%% file: posgap.tex
\documentclass[10pt]{article}

\usepackage{relsize}
\usepackage{tikz} 
\usepackage{url, amsmath,enumerate,fancyhdr,amssymb, amsthm, 
	pstricks, epsf, lscape, ifpdf, graphicx, floatrow, multirow} 
\usepackage{algorithm}
\usepackage{algorithmicx} 
\usepackage{algpseudocode}

\usepackage{import}
\usepackage{tabularx}
\newcolumntype{C}{>{\centering\arraybackslash}X}
\usepackage{amsmath}
\makeatletter
\newcommand{\leqnomode}{\tagsleft@true}
\newcommand{\reqnomode}{\tagsleft@false}
\makeatother

\usetikzlibrary{shapes.misc,shadows}

\usepackage[latin1]{inputenc}
\usepackage{enumerate} 
\usepackage[pagebackref=true, colorlinks=true, linkcolor=blue, citecolor=blue, linkcolor=blue, anchorcolor=red, urlcolor=blue]{hyperref}
\usepackage{url, amsmath,fancyhdr,amssymb, amsthm, 
	pstricks, epsf, lscape, ifpdf, graphicx, float}
\usepackage{algorithm}
\usepackage{color}
\usepackage{algpseudocode} 
\usepackage{tikz}
\usetikzlibrary{decorations.pathreplacing,angles,quotes}
\usetikzlibrary{matrix}
\usetikzlibrary{backgrounds, fit}

\algnewcommand\algorithmicinput{\textbf{Input:}}
\algnewcommand\INPUT{\item[\algorithmicinput]}
\algnewcommand\algorithmicinitialization{\textbf{Initialization:}}
\algnewcommand\INITIALIZATION{\item[\algorithmicinitialization]}

\usepackage{amsmath,amssymb,tikz-cd}

\usepackage{longtable,tabu}

\usepackage[margin=2.5cm]{geometry}

\usepackage{tabularx}

\renewcommand{\arraystretch}{1.2}

\newcommand{\pha}{\phantom}

\newcommand{\phao}{\phantom{0}}

\mathchardef\mhyphen="2D

\newcommand{\val}{\operatorname{val}}

\newtheorem{Definition}{Definition}

\newtheorem{Claim}{Claim}

\newtheorem{Example}{Example}
\newtheorem{Proposition}{Proposition}
\newtheorem{Lemma}{Lemma}
\newtheorem{Theorem}{Theorem}
\newtheorem{Corollary}{Corollary}
\newtheorem{Remark}{Remark}[section]
\newtheorem{Assumption}{Assumption}

\newcommand{\I}{\mathcal I}
\newcommand{\J}{\mathcal J}

\newcommand{\lin}{\operatorname{lin}}

\newcommand{\fr}{\operatorname{FR}}
\newcommand{\regfr}{\operatorname{REGFR}}

\newcommand{\eps}{\epsilon} 
\newcommand{\bpx}{\begin{pmatrix}}
	\newcommand{\epx}{\end{pmatrix}}
\newcommand{\bbx}{\begin{bmatrix}}
	\newcommand{\ebx}{\end{bmatrix}}

\newcommand{\bdef}{\begin{Definition}} 
	\newcommand{\commentout}[1]{}
	\newcommand{\co}[1]{}
	
	\newcommand{\nin}{\noindent}
	\newcommand{\ti}{\times}
	 
	\newcommand{\pf}[1]{\vspace{.35cm} \nin {\bf Proof {#1} }}

	\newcommand{\sym}[1]{{\cal S}^{#1}}
	\newcommand{\psd}[1]{{\cal S}_+^{#1}}
	\newcommand{\pd}[1]{{\cal S}_{++}^{#1}}
	\newcommand{\rad}[1]{\mathbb{R}^{#1}}

	\newcommand{\eref}[1]{(\ref{#1})}

	\newcommand{\la}{\langle}
	\newcommand{\ra}{\rangle}

	\newcommand{\beq}{\begin{equation}}
	\newcommand{\eeq}{\end{equation}}
	\newcommand{\beqa}{\begin{eqnarray}}
	\newcommand{\eeqa}{\end{eqnarray}}
	\newcommand{\ba}{\begin{array}}
		\newcommand{\ena}{\end{array}}
	\newcommand{\bac}{\begin{array}{ccccccccccc}}
		\newcommand{\eac}{\end{array}}
	\newcommand{\bprop}{\begin{Proposition}}
		\newcommand{\eprop}{\end{Proposition}}

	\newcommand{\beqast}{\begin{eqnarray*}}
		\newcommand{\eeqast}{\end{eqnarray*}}
	\newcommand{\benum}{\begin{enumerate}}
		\newcommand{\eenum}{\end{enumerate}}
	\newcommand{\bit}{\begin{itemize}}
		\newcommand{\eit}{\end{itemize}}
	\newcommand{\bth}{\begin{Theorem}}
		\newcommand{\enth}{\end{Theorem}}
	\newcommand{\ble}{\begin{Lemma}}
		\newcommand{\ele}{\end{Lemma}}
	\newcommand{\bex}{\begin{Example}}
		\newcommand{\eex}{\end{Example}}
	\newcommand{\bcor}{\begin{Corollary}}
		\newcommand{\ecor}{\end{Corollary}}
	\newcommand{\brem}{\begin{Remark}}
		\newcommand{\erem}{\end{Remark}}
	\newcommand{\bass}{\begin{Assumption}}
		\newcommand{\eass}{\end{Assumption}}

	\renewcommand{\arraystretch}{1.2}

	\newcommand{\bsmx}{\begin{small} \begin{pmatrix}}
			\newcommand{\esmx}{\end{pmatrix} \end{small}}
	
	\newcommand\bovermat[2]{%
		\makebox[0pt][l]{$\smash{\overbrace{\phantom{%
						\begin{matrix}#2\end{matrix}}}^{\text{#1}}}$}#2}
	
	\usepackage{mathtools}
	\makeatletter
	\renewcommand*\env@matrix[1][*\c@MaxMatrixCols c]{%
		\hskip -\arraycolsep
		\let\@ifnextchar\new@ifnextchar
		\array{#1}}
	\makeatother
	
	\usepackage{tikz,array,calc, xcolor}
	\usetikzlibrary{decorations.pathreplacing}
	
	\newcommand\tikzmark[2]{%
		\tikz[remember picture,overlay] 
		\node[inner sep=0pt,outer sep=2pt] (#1){#2};%
	}
	
	\newcommand\link[3]{%
		\begin{tikzpicture}[remember picture, overlay]
		\draw [decorate,decoration={brace,amplitude=5pt,raise=0pt}] 
		(#1.north west)--(#2.north east)node[above=3pt,midway]{#3};
		\end{tikzpicture}%
	}
	
	\title{\Large On Positive Duality Gaps in Semidefinite Programming} 
	\author{G\'{a}bor  Pataki} 

\begin{document}

		\maketitle

		\co{Figures:
			pic_simpleEx1v2.tex 
			pic_simpleEx23v4.tex
			pic_simpleEx4v5.tex
			psdwlabels.png
			} 
		
		\begin{abstract}
			
	We present a novel analysis of semidefinite programs (SDPs) with  positive duality gaps, 
	i.e. different optimal values in the primal and dual problems.
	These  SDPs  are extremely 
	pathological, often unsolvable, and  also serve as models of more general pathological convex programs. 
									However, despite their allure, they are not well understood even when they have 
								just two variables.
								
										We first completely characterize two variable SDPs with positive gaps; in particular,
		we transform them into a standard form 
		that  makes the positive gap  trivial to recognize.
				The transformation is very simple, as it mostly uses  
			elementary row operations 
			coming  from Gaussian elimination.
			We next show that the two variable case sheds light on larger SDPs with positive gaps: we present SDPs in any dimension in which the positive gap is 
			caused by the same structure as in the two variable case. We analyze a fundamental parameter,  
			the {\em singularity degree} of the duals of our SDPs,  and show that it is the largest that can 
			result  in a positive gap. 
			
			We finally generate a library of 
			difficult SDPs with positive gaps (some of these SDPs have only two variables) 
			and present a computational study.
			
		\end{abstract}


		
		{\em OR/MS subject classification:} Primary: convexity; secondary: programming-nonlinear-theory

\section{Introduction} \label{sect-intro}

In the last few decades we have seen an  intense growth of interest in semidefinite programs (SDPs), optimization problems with linear objective, linear constraints, and semidefinite matrix variables.
\co{
 Semidefinite programs (SDPs), optimization problems with linear constraints and objective, and positive semidefinite matrix variables,  are one of the most versatile and popular class of optimization problems of the last few decades. }

The recent appeal of SDPs is due to several reasons. First, SDPs 
are applied in areas as varied as combinatorial optimization, control theory, robotics,  polynomial optimization, and  machine learning.  Second, they naturally extend linear programming  (LP), 
and much research has been devoted to 
generalizing results from LP to SDPs, for example, to generalizing efficient interior point methods.
The extensive 
literature on SDPs includes textbooks, surveys and thousands of research papers. 

We consider an SDP in the form  
\begin{equation}\label{p}
\begin{split}
\sup & \,\, c^T x  \\
s.t.   & \,\,  \sum_{i=1}^m  x_i A_i \preceq B, \\
\end{split}\tag{\mbox{$P$}} 
\end{equation}
where  $A_1, \dots, A_m, \,$ and $B$ are $n \times n$ 
symmetric matrices,  $c \in \rad{m}$ is a vector, and for  symmetric matrices $S$ and $T$  we write 
$S \preceq T$ to say that $T - S$ is positive semidefinite (psd).  

To solve \eqref{p}, which we call the {\em primal} problem,  we rely on its natural {\em  dual} 
	\begin{equation}\label{d}
	\begin{split}
	\inf  & \,\, B \bullet Y  \\
	s.t. & \,\, A_i \bullet Y = c_i \, (i=1, \dots, m) \\
	& \,\, Y \succeq 0, 
	\end{split}\tag{\mbox{$D$}} 
	\end{equation}
	where the $\bullet$ product of symmetric matrices is the trace of their regular product.

SDPs inherit some of the duality theory of linear programs. For instance, if $x$ is feasible in \eqref{p} and $Y$ in \eqref{d},  then  the {\em weak duality inequality} 
$c^T x \leq B \bullet Y$
 holds. However, 
\eqref{p} and \eqref{d} 
 may not attain their optimal values,  and  their  optimal values  may even differ. In the latter case we say that there is a positive (duality) gap.

Among pathological SDPs, the ones with positive duality gaps may be the ``most pathological" or ``most interesting," depending on our point of view. 
They are in stark contrast with gapfree linear programs, and  
often look innocent, but still defeat SDP solvers.

 \bex \label{ex2} In the  SDP 
\beq \label{psmall} \tag{$P_{\rm small}$} 
\ba{rl}
\sup  &  x_2   \\
s.t. & x_1 \bpx 1 & 0 & 0 \\ 0 & 0 & 0 \\ 0 & 0 & 0 \epx + x_2 \bpx 0 & 0 & 1 \\ 0 & 1 & 0 \\ 1 & 0 & 0 \epx \preceq 
\bpx 1 & 0 & 0 \\ 0 & 1 & 0 \\ 0 & 0 & 0 \epx 
\ena
\eeq
the constraint is equivalent to 
\co{$$
\begin{pmatrix} 1 - x_1 & 0 & - x_2  \\
0         & 1 - x_2 & 0 \\
- x_2   & 0 & 0 \end{pmatrix} \succeq 0,
$$} 
 \begin{scriptsize} $\begin{pmatrix} 1 - x_1 & 0 & - x_2  \\
 	0         & 1 - x_2 & 0 \\
 	- x_2   & 0 & 0 \end{pmatrix} \succeq 0,$ \end{scriptsize} 
 \co{
 \begin{tiny} $\begin{pmatrix} 1 - x_1 & 0 & - x_2  \\
	0         & 1 - x_2 & 0 \\
	- x_2   & 0 & 0 \end{pmatrix} \succeq 0,$ \end{tiny}  } 
so 
 $x_2 =0$ always holds, and the optimal value of
 (\ref{psmall}) is $0.$ 

 Let $Y = (y_{ij}) \succeq 0$ be the dual variable matrix. By the first dual constraint 
 $y_{11}=0,$ so 
 the first row and column of $Y$ are  zero. Hence the dual is equivalent to 
\beq \label{problem-firstreduced}  
\ba{rrcl}  
\inf &   y_{22}   \\ 
s.t. &  y_{22}   & = & 1,  
\ena
\eeq
whose optimal solution is $1.$ 

This small example already shows that a  positive gap is bad news: the Mosek commercial SDP solver reports that 
(\ref{psmall}) is 
``primal infeasible." 

\eex
\begin{figure}[t] 
	\begin{center} 
		\subimport{imgs/}{pic_simpleEx1v2.tex}
	\end{center} 
	\label{figure-psmall}
	\caption{An illustration of Example \ref{ex2}} 
\end{figure}
We visualize this SDP in Figure \ref{figure-psmall}. The empty blocks in all matrices are zeroes (but a few zeroes are still shown to better visualize spacing). 
The nonzero  diagonal entries in all matrices are colored red.  In contrast, 
 we colored  the first row and column of $A_2$ blue. This blue portion of $A_2$ does not matter  when we compute $A_2 \bullet Y,$ where 
 $Y \succeq 0$ is feasible in the dual, since (as we just discussed) the first  row and column of $Y$ is zero. 
So the color scheme makes it clear that the dual is indeed equivalent to 
\eqref{problem-firstreduced}. 
 
The excitement about SDPs with positive gaps is evident  by the many examples published in surveys and textbooks: see for example, \cite{BonnShap:00, borwein2010convex, VanBo:96}. 
 
 Due to intensive research in the past few years, we now understand SDP pathologies much better, and can also remedy them, at least to some extent, and at least in theory. Among structural results, \cite{TunWolk:12} related positive gaps to complementarity in the homogeneous problems (with $B=0$ and $c=0$);  in  \cite{Pataki:17} we completely characterized pathological semidefinite {\em systems,} \co{
 with an unattained dual optimal value, or a positive duality gap 
 for {\em some}  $c \in \rad{m}.$ } and 
 \cite{lourencco2016structural, liu2019new} studied weakly  infeasible SDPs,  which are within zero distance of feasible ones. 
 
\co{
	
	3) I think that a citation in page 27 is not displayed correctly in the version you sent me.
	
	Best regards,
	Bruno
	}
 As to remedies, facial reduction algorithms \cite{BorWolk:81, Pataki:00B, WakiMura:12, Tuncel:11, Pataki:13, drusvyatskiy2017many,  lourencco2018facial} produce a dual problem which attains its optimal value and has zero gap with the primal. 
 Extended duals  \cite{Ramana:97,  KlepSchw:12, Pataki:13}  achieve the same goal and
  require no computation, but 
 involve extra variables (and constraints). 
 For the surprising connection of facial reduction and extended duals, see  
 \cite{RaTuWo:97, Pataki:00B, Pataki:13}.
 
 \co{Extended duals  \cite{Ramana:97,  KlepSchw:12, Pataki:13}  involve extra variables (and constraints), 
they attain their optimal value, and 
have zero duality gap with the primal. }

  In a broader context, a positive duality gap between \eqref{p} and \eqref{d} 
  implies zero distance to infeasibility, i.e., an arbitrarily small perturbation makes both infeasible \footnote{Of course, if any one of them is infeasible to start with, i.e., the duality gap is infinite, then 
  	this statement holds vacuously.}.  
  For literature on distance to  infeasibility see the  seminal paper  \cite{renegar1994some} and 
  many later works, e.g.,  \cite{freund2003complexity, pena2000understanding}.  
    Also note that SDPs are some of the simplest convex optimization problems with positive gaps, 
  so they serve as models 
  of other, more complex  pathological convex programs.
  An early prominent example 
  is Duffin's duality gap 
  \cite[Exercise 3.2.1]{borwein2010convex}; 
  \cite[Example 5.3.2]{bertsekas2009convex} is similar.

Despite the above studies and the plethora of published examples, positive duality gaps
do not seem well understood. 
For example, it is not difficult to see that  $m=1 \, $ implies no duality gap. 
However, even in the  $m=2$  case,  positive gaps have not yet been analyzed. 

\noindent{\bf Contributions} 
In a nutshell, 
we show 
that  simple certificates of positive gaps 
exist in a large class of SDPs, not just in artificial looking examples.

Our first result is 
\bth \label{thm-2var} 
Suppose $m=2.$ Then 
$\val(\ref{p}) < \val(\ref{d})$ iff 
$(\ref{p})$ has a reformulation
\beq \label{pref} \tag{\mbox{$P_{\rm ref}$}}
\ba{rllcl} 
\sup &    c_2' x_2  \\
s.t.   & 	 	x_1 \left(
\ba{c|c||c|c} 
\Lambda  & \phantom{\times} & \phantom{\times} &  \phantom{\times} \\ \hline
& & &  \\ \hline \hline 
& & &  \\ \hline
& & &   \\ 
\ena \right) & 
+ x_2 \left(
\ba{c|c||c|c} 
\ti & \ti & \ti &  M  \\ \hline
\ti & \Sigma  &  &   \\ \hline \hline 
\ti &  & - I_s &   \\ \hline
M^T & & &   \\ 
\ena \right) & \preceq & \left(\ba{c|c||c|c} 
I_p & \phantom{\ti} & \phantom{\ti} &  \phantom{\ti} \\ \hline
\phantom{\ti} & I_{r-p}  & \phantom{\ti} &  \phantom{\ti} \\  \hline  \hline  
\phantom{\ti} & \phantom{\ti} & \phantom{\ti} &  \phantom{\ti} \\ \hline 
\phantom{\ti} &  \phantom{\ti}  & \phantom{\ti} &  \phantom{\ti} \\ \hline 
\ena \right), 
\ena
\eeq
where $\Lambda$ and $\Sigma$ are diagonal, $\Lambda$ is positive definite, 
$M \neq 0, \,$ $c_2' > 0$ and $s \geq 0.$ 
\footnote{Example \ref{ex2} needs no reformulation and has $\Sigma = [1]$ and $s=0.$}.  
\enth
\qed

Hereafter,  $\val()$ denotes   
the optimal value of an optimization problem. 
We partitioned the matrices 
to show  their order, e.g., $\Lambda$ has order $p.$ 
The empty blocks are zero, and the ``$\times$" blocks may  have arbitrary elements.

In Subsection
\ref{subsection-preliminaries} 
we precisely define ``reformulations."  
However, if we believe that 
a reformulated problem 
has a positive gap with its dual iff the original one does, we can already prove the ``easy",  the ``If" direction of Theorem \ref{thm-2var}. 
To build intuition, we give this proof below; it  
essentially reuses the argument from Example \ref{ex2}.

\pf{of ``If" in Theorem \ref{thm-2var}:} 
Since $M \neq 0, \,$ we have $x_2 = 0$ in any feasible solution of
\eref{pref} so   $\val(\ref{pref})=0.$ 

Suppose  next that $Y$ is feasible in the dual of \eref{pref}. By the first dual constraint 
 a positively weighted  linear combination of the first $p$ diagonal elements of $Y$ is zero, so these elements are all zero. Since $Y \succeq 0, \, $ the first $p$ rows and columns of $Y$ are also zero.
 We infer that the dual is equivalent to the reduced dual 
\beq \label{dred} \tag{$D_{\rm red}$} 
\ba{rrcl}
\inf & \bpx I & 0 \\ 0 & 0 \epx \bullet Y' \\
s.t. & \bpx \Sigma & 0 \\ 0   & - I_s  \epx \bullet Y'  & =            & c_2'  \\
&                                                                  Y'  & \succeq & 0.
\ena
\eeq
If $\Sigma \leq 0,$  then \eqref{dred} is infeasible, so its optimal value is 
$+ \infty.$ 
Otherwise, some diagonal element of $Y'$ corresponding to a diagonal element of $\Sigma$ 
must be positive, so 
$\val \eqref{dred}$ is positive and finite. 

Either way, there is a positive duality gap.
\qed

After reviewing preliminaries in subsection \ref{subsection-preliminaries}, 
in 	Section \ref{section-2var} we prove the ``only if" direction of  Theorem \ref{thm-2var} and some natural corollaries.
    For example,   Corollary \ref{corollary-2var} shows 
 that when $m=2,$ 
  the ``worst" pathology -- positive gap coupled with unattained primal or dual optimal value -- 
  is entirely absent. 
  
  We next show that the two variable case,  although it may seem much too special, 
    sheds light on larger SDPs with positive gaps. 
\co{We next show that, although  the two variable case may seem much too special, it 
sheds light on larger SDPs with positive gaps.} 
  Section \ref{section-similar-in-anydimension} presents SDPs in any dimension,  in which the positive gap is 
  manifested  by the same structure as in the two variable case. 
    In Section \ref{section-singdegree}  we compute 
  the  singularity degree  -- the minimum number of steps that a facial reduction algorithm needs to regularize an SDP  -- 
  of the duals of our SDPs.
  The first dual is just 
  (\ref{d}) and the second is the  homogeneous dual 
  \beq \label{HD} \tag{$HD$}
  \begin{array}{rcl} 
  	A_i \bullet  Y & = & 0 \,\,  (i=1, \dots, m)  \\ 
  	B   \bullet Y   & = & 0  \\
  	Y & \succeq &  0.
  \end{array} 
  \eeq
  We show that the singularity degrees of \eqref{d} and \eqref{HD} corrresponding to the SDPs in Section 
  \ref{section-similar-in-anydimension} are  $m-1$ and $m, \,$ respectively \footnote{Precisely, in the ``single sequence" SDPs in Section \ref{section-similar-in-anydimension}, the singularity degree of 
  	\eqref{d} is  $m-1.$ In  the ``double sequence"  SDPs, the singularity degree of 
  	\eqref{d} is  $m-1,$ {\em and}  the singularity degree of \eqref{HD} is $m.$ }. 
  Section \ref{section-maximal-deg} is a counterpoint: 
  it shows that 
   the singularity degrees of (\ref{d})   and of (\ref{HD})  are always
  $\leq m$ and $\leq m+1, \,$ respectively, and when equality holds there is no duality gap. 
  
  Finally, in Section \ref{section-computational} we  generate a library of SDPs, in which  the positive gap can be verified 
  by simple inspection and in exact, integer  arithmetic. However, for  current software 
  these SDPs turn out to be essentially unsolvable.

We believe that some of the paper's results are of independent interest. 
For example, Theorem \ref{theorem-d-equals-mplus1} shows that {\em maximal} singularity degree (equal to $m$) in \eqref{HD} implies {\em minimal} singularity degree (equal to $0$) in 
\eqref{d}. 
We also expect the visualizations of Figures \ref{figure-psmall}, \ref{figure-sparsity-E_ij} and 
\ref{figure-double-sequence} to be useful  to others. 

\noindent{\bf Related literature} 
By Theorem 5.7 in  \cite{TunWolk:12}  when the ranks of the maximum rank solutions sum to $n-1 \, $ in the 
homogeneous primal-dual pair (with $c=0$ and $B=0$),  
 \eqref{p} and \eqref{d} have a positive duality gap for a suitable $c$ and $B.$ 

In  \cite{Pataki:17}  we characterized pathological semidefinite {\em  systems}, which have an unattained dual value or positive gap for {\em some}  $c \in \rad{m}.$ 
However, \cite{Pataki:17} cannot distinguish among 
``bad" objective functions. For example, 
it cannot tell which $c$  gives a positive gap, and which gives zero gap and unattained dual value, 
a much more harmless pathology. 
 
Weak infeasibility of \eqref{d} (or by symmetry, of \eqref{p}) 
is the same as  an infinite duality gap  in all salient cases, as 
we show in Proposition \ref{prop-weakinfeas}. See  \cite{lourencco2016structural} 
for a proof that any weakly infeasible SDP contains a ``small"  such SDP of dimension at most $n-1.$ 
\co{In fact, in Section \ref{section-maximal-deg} we use some ideas from 
 \cite{lourencco2016structural}. } 
Furthermore, 
\cite{LiuPataki:15} and \cite{liu2017exact} characterized infeasibility and weak infeasibility in conic LPs, by {\em reformulating} them into standard forms. We will use the same technique in this work. 

The minimum number of steps that a facial reduction algorithm needs is 
the  {\em singularity degree} of the SDP, a  fundamental parameter introduced in   \cite{Sturm:00} and used 
in many later works. 
	 				An upper bound on the singularity degree of an  SDP with order $n$ matrices is $n-1, \,$ and 
			this bound is tight, as a nice  example  in  \cite{Tuncel:11} shows.  
				On the other hand, the SDPs in this paper are the first ones that have such a large singularity degree, 
				a positive duality gap and a structure  inherited from the two variable case.  
					We further refer to \cite{lourencco2018facial} for an improved bound on the singularity degree, 
					when  the underlying cone has polyhedral faces;   to  \cite{lourencco2017amenable} for a broad generalization of the error bound  of \cite{Sturm:00} to conic LPs over {\em amenable cones};
						and to \cite{permenter2014partial} and \cite{zhu2019sieve} for recent implementations of facial reduction. 
		 \co{Recent implementations of facial reduction are Some recent implementations of facial reduction algorithms are given in \cite{permenter2014partial} 
		 \cite{zhu2019sieve}} 
		 
		  
		 \co{To implement a facial reduction algorithmm as described in 
		 one must solve SDPs in exact arithmetic. However, recent ``partial" im
		 The theory of facial reduction is well established, however,  these algorithms  
		 are not easy to implement. For recent 	} 
In other related work,  \cite{de2000self}  used  self-dual embeddings and Ramana's dual to recognize SDP pathologies. Their method has not been implemented yet, 
as it must solve SDP subproblems in exact arithmetic. 

\noindent{\bf Reader's guide} Most of the paper (in particular, all of 
Sections \ref{section-2var}, \ref{section-similar-in-anydimension} 
and \ref{section-computational})  can be read with a minimal background in linear algebra and semidefinite 
programming, all of which we summarize in Subsection \ref{subsection-preliminaries}. 
The proofs are  short and fairly elementary, they mostly use only elementary linear algebra, 
and we  illustrate our results with many examples.

\subsection{Preliminaries} 
\label{subsection-preliminaries} 
\paragraph{Reformulations} 
We first introduce reformulations, a tool that we used in recent works \cite{Pataki:17, LiuPataki:15} to analyze pathologies in SDPs.  
\begin{Definition} \label{definition-reform}  We say that we  {\em reformulate } the pair of SDPs (\ref{p}) and (\ref{d})  if we apply to them some of the following operations:
	\benum
	\item \label{slack} Replace $B$ by $B +\lambda A_j, \,$ for some $j$ and $\lambda \neq 0.$ 
	\item \label{exch} Exchange $(A_i, c_i)$ and $(A_j, c_j), \,$ where $ i \neq j.$ 
	\item \label{trans} Replace $(A_i, c_i)$ by $\lambda (A_i, c_i) + \mu (A_j, c_j), \,$ where $\lambda \neq 0.$
	\item \label{rotate} Apply a similarity transformation 
	$T^T()T$ to all $A_i$ and $B, \,$ where $T$ is an invertible matrix.  
	\eenum
	We also say that by reformulating (\ref{p}) and (\ref{d}) 
	 we obtain a {\em reformulation}. 
\end{Definition}
(Of course, we can apply any of these operations, in any order.) 

Operations $(\ref{slack})\mhyphen (\ref{trans})$  correspond to elementary row operations 
done on  (\ref{d}). 
For example,  operation (\ref{exch}) exchanges the constraints 
$$
A_i \bullet Y = c_i \,\, {\rm and} \,\, A_j \bullet Y = c_j.
$$
Clearly, (\ref{p}) and (\ref{d}) attain their optimal values iff they do so after a reformulation, and 
reformulating them  also preserves duality gaps.

\paragraph{Matrices} 
As usual, $\sym{n}, \psd{n},$ and $\pd{n}$ stand for the set of $n \times n$ symmetric, symmetric positive semidefinite (psd), and symmetric positive definite (pd) matrices, respectively. For $S, T \in \sym{n}$ we write $T \succ S$ to say $T - S \in \pd{n}.$  

We denote by $I_r$ and $0_r$ the $r \times r$ identity matrix and all zero matrix, respectively.
For matrices $A$ and $B$ we denote their concatenation along the diagonal by $A \oplus B, \, $ i.e., 
$$
A \oplus B := \bpx A & 0 \\
0  &  B \epx.
$$
Accordingly,  $\psd{r} \oplus 0_s$ is  the set of order $r+s$ symmetric matrices whose 
upper left $r \times r$ block is psd, and the rest is zero. Sometimes we just write $\psd{r} \oplus 0$ 
if the order of  the zero part is clear from the context. 
The meanings of $\pd{r} \oplus 0_s$ and $0_s \oplus \psd{r}$ are similar. 

\paragraph{Strict feasibility (or not)} 
We call $Z \succeq 0$ a {\em slack matrix in (\ref{p})}, if 
$Z = B - \sum_{i=1}^m x_i A_i$ for some $x \in \rad{m}.$ 

We say that (\ref{p}) is strictly feasible \footnote{or it satisfies Slater's condition},  
if there is 
a positive definite slack in it. If (\ref{p}) is strictly feasible, then there is no duality gap, 
and $\val(\ref{d})$ is attained when finite.
Similarly, we say that (\ref{d}) is strictly feasible, if    
it has a positive definite feasible $Y.$ In that case there 
is no duality gap, and $\val(P)$ is attained when finite. 

Can we certify the {\em lack} of strict feasibility? Given an affine subspace $H \subseteq \sym{n}$ such that  $H \cap \psd{n}$ is nonempty, 
the Gordan-Stiemke theorem for the semidefinite cone states 
\beq \label{eqn-HKlemma} 
H\cap \pd{n}  =\emptyset\Leftrightarrow H^{\perp}\cap ( \psd{n} \setminus \{ 0 \})  \neq \emptyset.
\eeq
For example, if  $H = \{ Z: \, Z = B - \sum_{i=1}^m x_i A_i  \, {\rm for \, some \,}   x \in \rad{m} \, \}, $ 
then $H \cap \psd{n}  \, $ is the set of feasible slacks in \eqref{p}. So \eqref{p} is 
not strictly feasible, iff its homogeneous dual \eqref{HD} has a nonzero solution. 
 
 We make  the following
\bass \label{assumption-slack} 
Problem (\ref{p}) is feasible, the $A_i$ and $B$ are linearly independent, 
$B$ is of the form 
\begin{equation} \label{equation-Bslack} 
B = \bpx I_r & 0 \\ 0 & 0 \epx, \,\, {\rm where} \,\,  0 \leq r < n,
\end{equation}
and it is the maximum rank slack in $(P).$ 
\eass
The assumption about $B$ is  easy to satisfy, at least in theory. The argument goes as follows:  
suppose $Z$ is a maximum rank slack 
in (\ref{p}) and $Q$ is a matrix of  suitably scaled eigenvectors of $Z.$ If we first replace
$B$ by $Z \,$  then replace 
$A_i$ by $Q^T A_i Q$ for all $i, \,$ and $B$ by $T^T B T, \,$ then $B$ will be in the required form.

\paragraph{Schur complement condition for positive (semi)definiteness} 
We recap a classic condition for positive (semi)definiteness. If $G \in \sym{n}$ is partitioned as  
$$
G = \bpx G_{11} & G_{12} \\
              G_{12}^T & G_{22} \epx
$$
with $G_{22} \succ 0, \, $ then the following equivalences hold: 
\begin{eqnarray} \label{schur-psd} 
G \succeq 0 &  \Leftrightarrow & G_{11} - G_{12} G_{22}^{-1} G_{12}^T \succeq 0,  \\ \label{schur-pd}
G \succ 0 &  \Leftrightarrow & G_{11} - G_{12} G_{22}^{-1} G_{12}^T \succ 0. 
\end{eqnarray}

\section{The two variable case} 
\label{section-2var}

\subsection{Proof of ``Only if" in Theorem \ref{thm-2var}} 
	
We now turn to the proof of the ``Only if" direction in Theorem 
\ref{thm-2var}. We chose a proof that employs the minimum amount of convex analysis, 
namely only the Gordan-Stiemke theorem  \eqref{eqn-HKlemma}. The rest of the proof is just linear algebra.

The main idea 
is that (\ref{d}) cannot be strictly feasible, 
otherwise the duality gap would be zero. We first make the lack of 
strict feasibility obvious by creating the constraint 
\begin{equation} \label{eqn-lambda} 
(\Lambda \oplus 0) \bullet Y = 0,
\end{equation}
where 
$\Lambda$ is diagonal, with  positive diagonal entries.   
Clearly, if $\Lambda$ is $p \times p \,$ and $Y \succeq 0$ satisfies \eqref{eqn-lambda}, then the first $p$ rows and columns of $Y$ must be zero.

To create the constraint \eqref{eqn-lambda}, we first perform  a  
 facial reduction step (using the Gordan-Stiemke theorem   \eqref{eqn-HKlemma}), 
 then a reformulation step. 
\co{Since we only need one facial reduction step, 
we simply invoke  
the Gordan-Stiemke theorem once. } 
We next analyse cases to show that the second constraint matrix must be in a certain form, 
 and further reformulate (\ref{p}) to put it into the final form \eqref{pref}. 

We need a basic lemma, whose proof is in Appendix \ref{appendix-proof-lemma-rotate}.

\ble \label{lemma-rotate} 

Let 
$$
G = \bpx G_{11} & G_{12} \\ G_{12}^T & G_{22} \epx,
$$
where	$G_{11} \in \sym{r_1}, \, G_{22} \in \psd{r_2}.$ 

Then there is an invertible matrix $T$ such that 
$$
T^TGT = \left( \ba{c|cc} \Sigma & 0     & W \\ \hline 
0          & I_s   & 0 \\
W^T     & 0     & 0 \ena \right) \; {\rm and} \; T^T \bpx I_{r_1} & 0 \\ 0 & 0 \epx T = \bpx I_{r_1} & 0 \\ 0 & 0 \epx, 
$$
where $\Sigma \in \sym{r_1}$ is diagonal and  $s \geq 0.$ 

\ele
\qed

\pf{of "Only if" in Theorem \ref{thm-2var}} 
\co{Throughout the reformulation process, we call the 
primal and dual problems (\ref{p}) and (\ref{d}), and we  also call the constraint matrices on the
left $A_1'$ and $A_2'.$ We start with $A_1' = A_1$ and $A_2' = A_2.$ 
} 
We call the 
primal and dual problems (\ref{p}) and (\ref{d}), and  the constraint matrices on the
left $A_1'$ and $A_2'$ throughout the reformulation process. We start with $A_1' = A_1$ and $A_2' = A_2.$  

\paragraph{Case 1: (\ref{d}) is feasible} 
.

We break the proof into four parts: facial reduction step and first reformulation; transforming 
$A_1'; \,$ transforming $A_2';  \,$ and ensuring $c_2' > 0.$ 

\paragraph{Facial reduction step and first reformulation} 
Let 
\beq \label{eqn-defineH} 
\begin{array}{rcl}
	H & = & \{ \, Y \, | \, A_i \bullet Y = c_i \, \forall i \, \} \\
	& = &  \{ \, Y \, | \, A_i \bullet Y = 0 \, \forall i \, \} + Y_0,   
\end{array}
\eeq
where $Y_0 \in H$ is arbitrary. Then the feasible set of the dual \eqref{d} is $H \cap \psd{n}.$ 
Since (\ref{d}) is not strictly feasible,   by the Gordan-Stiemke theorem (\ref{eqn-HKlemma}) there is 
$$
\ba{rcl}
A_1' & \in & \bigl( \psd{n} \setminus \{ 0 \} \bigr) \cap H^\perp. 
\ena
$$
Thus   for some $\lambda_1$ and $\lambda_2$ reals we have 
$$
\ba{rcl} 
 A_1' & = & \lambda_1 A_1 + \lambda_2 A_2,  \, {\rm and}  \\ 
A_1' \bullet Y_0 & = & 0. 
\ena
$$
Since 
$$
\lambda_1 c_1 + \lambda_2 c_2 \, = \,  (\lambda_1 A_1 + \lambda_2 A_2) \bullet Y_0 \, = \,  A_1' \bullet Y_0 \, = \, 0,
$$
we can reformulate the feasible set of (\ref{d})\footnote{I.e.,  we just ignore $B.$}  using only operations (\ref{exch}) and (\ref{trans}) in Definition 
\ref{definition-reform} as 
\beq \label{eqn-Dref-1st} 
\ba{rcl}
A_1' \bullet Y & = & 0 \\
A_2' \bullet Y & = & c_2'  \\
Y & \succeq & 0
\ena
\eeq
with some $A_2'$ matrix and $c_2'$ real number. (To be precise, 
if $\lambda_1 \neq 0$ then we multiply the first  equation in \eqref{d} by $\lambda_1$ and add $\lambda_2$ times the second equation to it. If $\lambda_1 = 0,$ then $A_1' \neq 0$ implies that $\lambda_2 \neq 0$ and $c_2 = 0, \, $ so  we multiply the second equation in \eqref{d} by $\lambda_2$ and exchange it with the first. )

\paragraph{Transforming $A_1'$} 
Since $A_1' \succeq 0$ and $B$ is the maximum rank slack in $(P),$ 
the only nonzero entries of $A_1'$ are in its upper left
$ r \times r$ block,  otherwise $B - x_1 A_1'$ would be a slack with larger rank than 
$r$ for $x_1 <0.$ 

Let $p$ be the rank of $A_1', \,$  $Q$  a matrix of length 1 eigenvectors of the upper left $r \times r$ block of $A_1', $ set $T = Q \oplus I_{n-r},$ and apply the transformation 
$T^T()T$ to  $A_1', A_2' \,$ and $B.$ After this $A_1'$ looks like 
\beq \label{equation-A1prime} 
A_1'  \, = \, \left( 
\ba{c|c||cc} 
\Lambda & \phantom{\times} & \phantom{\times} &  \phantom{\times} \\ \hline
& & &  \\ \hline \hline 
& & &  \\ 
& & &   \\ 
\ena \right),
\eeq
where $\Lambda \in \sym{p}$ is diagonal with positive diagonal entries. 
From now the upper left $r \times r$ corner of all matrices will 
be bordered by  double lines. 
Note that $B$ is still in the same form as in the beginning (see Assumption  \ref{assumption-slack}). 

\paragraph{Transforming $A_2'$} 
Let $S$ be the lower $(n-r) \ti (n-r)$ block of $A_2'.$ We claim that 
\begin{equation} \label{eqn-S-not-indefinite} 
S \, {\rm cannot \, be \, indefinite,}
\end{equation}
so suppose it is.
Then the equation $S \bullet Y' = c_2'$ has a positive definite solution $Y'.$ Then 
\beq \nonumber 
Y \, := \, \left( 
\ba{c||c}
0 & 0 \\ \hline \hline 
0 & Y' 
\ena
\right)
\eeq
is feasible in \eqref{d}  with value $0, \,$ thus
$$
0 \, \leq \, \val \eqref{p}  \, \leq \, \val(\ref{d}) \, \leq \, 0,
$$
so the duality gap is zero, which is a contradiction.  We thus proved \eqref{eqn-S-not-indefinite}. 

We can now assume $S \succeq 0$ (if $S \preceq 0,$ we just multiply $A_2'$ and $c_2'$ by $-1 \,$). 
Recall that $\Lambda$ in (\ref{equation-A1prime}) is 
$p \times p, \,$ where $p \leq r.$ Next we apply Lemma \ref{lemma-rotate} with 
\begin{eqnarray*}
	G & := & {\rm lower \, right \,} (n-p) \times (n-p) \, {\rm block \, of \,} A_2', \\
	r_1 & := & r - p, \,\, {\rm and} \\
	r_2 & := & n -r. 
\end{eqnarray*}
Let $T$ be the invertible matrix supplied by Lemma \ref{lemma-rotate}, and apply the 
transformation  \mbox{$(I_p \oplus T)^T() (I_p \oplus T)$} to $A_1', A_2'$ and $B.$ 
This operation keeps $A_1'$ as it was.
It also keeps $B$ as it was, since the transformation  $T^T()T$ keeps $(I_{r-p} \oplus 0)$ the same. 

Next we  multiply both $A_2'$ and $c_2'$ by $-1$ to make $A_2'$ look  like 
\beq \label{equation-finalA2} 
A_2' \, = \, \left(
\ba{c|c||c|c} 
\ti & \ti & \ti &  M  \\ \hline
\ti & \Sigma & 0 &   W \\ \hline \hline 
\ti & 0 & - I_s &   \\ \hline
M^T & W^T & &   \\ 
\ena \right) \; {\rm for \; some \;} M \; {\rm and} \; W. 
\eeq
We next claim 
\begin{equation} \label{eqn-Mneq0-or-Wneq0} 
W \neq 0 \,\, {\rm or} \,\, M \neq 0.  
\end{equation}
For the sake of obtaining a contradiction, 
suppose that $W = 0$ and $M = 0.$ Let $(P_{\rm red})$ be the SDP obtained from \eqref{p} 
by deleting the first $p$ rows and columns from all matrices. Since $\Sigma$ is diagonal, 
$(P_{\rm red})$ is just a linear program with one constraint. Hence its set of feasible solutions is some closed interval,  say $[\alpha, \beta], \,$ 
and it has an optimal solution $x_2 = \alpha$ or $x_2 = \beta.$ 
 Further, its  dual is equivalent to \eqref{dred}.
  
  We next claim 
 \begin{equation} \label{eqn-valp=vald} 
  \val \eqref{p} \, = \, \val P_{\rm red})  \, = \, \ \val \eqref{dred} \, = \, \val \eqref{d}.
 \end{equation}
 Indeed, the first equation in \eqref{eqn-valp=vald} follows since for any $x_2 \in (\alpha, \beta)$ 
 the matrix 
 $$
 \bpx I - x_2 \Sigma & 0 \\
 0                                 & x_2 I_s 
 \epx
 $$
 is positive definite. Hence for any such $x_2$  the upper left order $r+s$ block of $Z := B - x_1 A_1 - x_2 A_2$ is positive definite if  
 $x_1$ is sufficiently negative: this follows from the Schur-complement condition for positive definiteness
 \eqref{schur-pd}. 
 
 The second equation in \eqref{eqn-valp=vald} follows since $(P_{\rm red})$ and \eqref{dred} are 
 linear programs; and the third equation follows since in any $Y$  feasible in \eqref{d} the first $p$ rows and columns are zero. 

 In summary, from the assumption that $M = 0$ and $W = 0$ we deduced that the duality gap is zero, 
 which is a contradiction. This 
 proves  \eqref{eqn-Mneq0-or-Wneq0}. 

We next claim that 
\begin{equation} \label{eqn-x2=0} 
x_2 = 0 \,\,  {\rm  in \, any \, feasible \, solution \, of \, } \eqref{p}. 
\end{equation}
Indeed, suppose $(x_1, x_2)$ is feasible in \eqref{p}, 
and  $x_2 \neq 0.$ Since $M \neq 0$ or $W \neq 0, \, $ the corresponding slack matrix has a $0$ diagonal entry, and a corresponding nonzero offdiagonal entry, thus it cannot be psd, which is a contradiction. 
We thus proved \eqref{eqn-x2=0}. 

Since $x_2 = 0$ always holds in \eqref{p}, we deduce $\val \eqref{p} =0.$ 

Next we claim
\begin{equation} \label{eqn-W=0}
W = 0, 
\end{equation}
so suppose  $W \neq 0.$ Then 
we define 
$$
Y \, = \, \left(
\ba{c|c||c|c} 
\phantom{\ti} & \phantom{\ti} & \phantom{\ti} &  \phantom{M}   \\ \hline
\phantom{\ti} & \eps I  & \phantom{V} &   * \\ \hline \hline 
\phantom{\ti} & \phantom{V^T}  &\phantom{ I_s} &   \\ \hline
\phantom{M^T} & * & &  \lambda  I  \\ 
\ena \right),
$$
where $\eps> 0, \,$ we choose the  ``*" block so that $A_2' \bullet Y = c_2', \,$ 
we choose $\lambda >0$ large enough to ensure $Y \succeq 0,$ and 
the empty blocks of $Y$ as  zero.
Consequently, $B \bullet Y = (r-p) \eps, \,$ so letting $\epsilon  \searrow 0$ we deduce $\val(\ref{d})=0, \,$ which is a contradiction.
This proves \eqref{eqn-W=0}.

\paragraph{Ensuring $c_2' > 0.$}  We  have  $c_2' \neq 0, \,$ otherwise the primal objective function would be $(0,0), \, $ so 
the duality gap would be zero. 

First, suppose $s > 0.\,$ 
 We will prove that in this case $c_2' > 0$ must hold, so to obtain a contradiction, assume $c_2' < 0.$ Let 
$$
Y \, := \, \left(
\ba{c|c||c|c} 
\phantom{\ti} & \phantom{\ti} & \phantom{\ti} &  \phantom{M}   \\ \hline
\phantom{\ti} & \phantom{\Lambda}   & \phantom{V} &   \phantom{W}  \\ \hline \hline 
\phantom{\ti} & \phantom{V^T}  & (-c_2'/s)  I_s  &   \\ \hline
\phantom{M^T} & \phantom{W^T}  & &  \phantom{\ti}   \\ 
\ena \right),
$$
where, as usual, the empty blocks are zero. Then $Y$  is feasible in (\ref{d}) with value $0, \,$ which is a  contradiction, and proves $c_2' > 0.$ 

Next, suppose $s=0.$ If $c_2' > 0, \,$ then we are done; if 
$c_2' < 0, \,$ then we multiply both $A_2'$ and $c_2'$ by $-1$ to ensure $c_2' > 0.$ 

We have thus transformed \eqref{p} into the form of \eqref{pref} and this completes the proof of Case 1. 

\paragraph{Case 2: (\ref{d}) is infeasible}  Since there is a 
positive duality gap, we see that  $\val(P) < + \infty.$ 

Consider the SDP 
\beq \label{eqn-Ylambda} 
\ba{rrcl} 
\inf & - \lambda \\
s.t.   & A_i \bullet Y - \lambda c_i & = & 0 \,\,  \forall i \\
&                                     Y & \succeq & 0 \\
&                          \lambda & \in & \rad{}, 
\ena
\eeq 
whose optimal value is zero: indeed, if $(Y, \lambda)$ were feasible in it with $\lambda > 0,$ then
$(1/\lambda)Y$ would  be  feasible in \eqref{d}. 
We next claim that 
\beq \label{eqn_Y-DNE} 
\not \exists (Y, \lambda) \, {\rm feasible \, in \, }  (\ref{eqn-Ylambda}) \,  {\rm \, such \, that \, } Y \succ 0 ,
\eeq
so suppose there is such a $(Y, \lambda).$ 
We next construct an  SDP  in the standard dual form, which is equivalent to \eqref{eqn-Ylambda}:
\beq \label{eqn-Ylambda-2}  
\ba{rrcl} 
\inf  & \bpx 0_n  & & \\  & -1 & \\ & & 1 \epx  \bullet \bar{Y} \\
s.t.  & \bpx A_i & & \\  & - c_i & \\  & & c_i \epx  \bullet \bar{Y} & = & 0 \,\,  \forall i \\
& \bar{Y}    & \in & \psd{n+2}.
\ena
\eeq
Observe that the $\lambda$ free variable in \eqref{eqn-Ylambda}  is split as 
$
\lambda \,  = \,  \bar{y}_{n+1, n+1} - \bar{y}_{n+2, n+2} 
$
in \eqref{eqn-Ylambda-2}, where $\bar{y}_{ii}$ is  the  $(i,i)$th diagonal elements of $\bar{Y}.$ 

Thus, \eqref{eqn-Ylambda-2} is strictly feasible with $\bar{Y} = Y \oplus  ( \bar{y}_{n+1, n+1} ) \oplus ( \bar{y}_{n+2, n+2}), \, $ 
 where 
$(Y, \lambda)$ is feasible in \eqref{eqn-Ylambda} with $Y \succ 0, \,$ and 
$\bar{y}_{n+1, n+1}$ and $\bar{y}_{n+2, n+2}$  are positive reals  whose difference is $\lambda.$ 

In summary, 
\eqref{eqn-Ylambda-2} is strictly feasible, and 
has zero optimal value (because \eqref{eqn-Ylambda} does). 
So the dual of \eqref{eqn-Ylambda-2} is feasible, i.e., there is $\bar{x} \in \rad{m}$ s.t. 
\beq \label{eqn-xbar} 
\ba{rcl} 
\sum_{i=1}^m \bar{x}_i A_i & \preceq & 0 \\
\sum_{i=1}^m \bar{x}_i c_i & = & 1.
\ena  
\eeq
Adding a  large multiple of $\bar{x}$ to a feasible solution of (\ref{p}) we deduce 
$\val(\ref{p}) = + \infty, \,$ which is a contradiction.
We thus proved  (\ref{eqn_Y-DNE}).

Of course, (\ref{eqn_Y-DNE}) means $\lin H \cap \pd{n} = \emptyset,$ where 
$H$ is defined in (\ref{eqn-defineH}). 
Since $(\lin H)^\perp = H^\perp, \,$ we next invoke the Gordan-Stiemke theorem 
(\ref{eqn-HKlemma}) with $\lin H$ in place of $H$ and 
complete the proof just like we did in Case 1. 
\qed 

\co{
\brem {\rm Among the several possible proofs of Theorem \ref{thm-2var} we chose one that employs the minimum amount of convex analysis, 
namely only the Gordan-Stiemke theorem  \eqref{eqn-HKlemma}. The rest of the proof is just linear algebra.
\co{, and 
we hope that our treatment  makes  the proof accessible to the broadest 
possible audience. } 

We could shorten the proof by invoking other results in SDP duality. For example Theorem 1     in ??? bad SDP short shows that 
there is a $V$ linear combination of $A_1'$  and $A_2'$ of the form 
\begin{equation} \label{eqn-V} 
V = \bpx V_{11} & V_{12} \\
       V_{12}^T & V_{22} \epx, 
\end{equation}
where $V_{22} \in \psd{n-r}$ and the rangespace of $ V_{12}^T$is not contained in the rangespace of $V_{22}.$ 
	
	Once we are done with Step ``{\bf Transforming $A_1'$,} the matrix $A_2'$ itself  must be of the form as $V$ is in \eqref{eqn-V}, and this allows us to make some shortcuts. 
	}
\erem
} 
\subsection{Some corollaries} 
\label{subsection-corollaries} 

Arguably the worst  possible pathology of SDPs is a positve duality gap accompanied by an unattained 
primal or dual optimal value. 
Luckily, as we next show, 
this worst pathology 
does not happen when $m=2.$ 

\begin{Corollary} \label{corollary-2var} 
Suppose $m=2, \,$  (\ref{p}) is feasible, and  $\val(\ref{p}) < \val(\ref{d}).$ 
Then (\ref{p}) attains its optimal value, and so does (\ref{d}) if it is feasible. 
\end{Corollary}
\nin {\bf Proof} Assume  the conditions above hold and assume w.l.o.g. 
that we reformulated (\ref{p})  into (\ref{pref}) and \eqref{d} into $(D_{\rm ref})$, the dual of (\ref{pref}). 
We will prove the above statements for \eqref{pref} and $(D_{\rm ref}).$

Since $x_2 = 0$ always holds in (\ref{pref}), its optimal value is $0$ and it is attained. 

Assume that $(D_{\rm ref})$ is feasible. 
Then it  is equivalent to the reduced dual
\eqref{dred} in which 
the matrix  $\Sigma$ is diagonal. So (\ref{dred})  is just a linear program, which attains its optimal value, hence so does  
$(D_{\rm ref}).$
\qed

	\brem{\rm  When the assumptions of Corollary \ref{corollary-2var} hold, we can prove 
		an even stronger result. In that case 
	 the objective value of \eqref{pref} is identically zero over the feasible set.  Since we obtained \eqref{pref} from \eqref{p} using only operations 
		\eqref{exch}$\mhyphen$\eqref{rotate}, the same holds for \eqref{p}.
	}
	\erem

We now turn to studying the 
semidefinite {\em system} 
\beq \label{p-sd} \tag{\mbox{$P_{SD}$}}
\sum_{i=1}^m x_i A_i \preceq B.  
\eeq
In \cite{Pataki:17} we characterized when (\ref{p-sd}) is {\em badly behaved}, 
meaning when there is 
$c \in  \rad{m}$ such that (\ref{p}) has a finite optimal value, but (\ref{d}) has 
no solution with the same value. Hence we may wonder, when is there 
$c \in  \rad{m}$ that leads to a positive gap, i.e., when  is (\ref{p-sd})  
``really" badly behaved?

The following straightforward corollary of Theorem \ref{thm-2var} 
settles this question when $m=2.$ It relies on  {\em reformulating} 
(\ref{p-sd}), i.e., reformulating (\ref{p}) with {\em some} $c.$ 

\bcor  \label{corollary-2var-system} 
Suppose $m=2.$ Then 
there is $(c_1, c_2)$ such that 
$\val(\ref{p}) < \val(\ref{d})$ iff 
(\ref{p-sd}) has a reformulation
\beq \nonumber 
\ba{rllcl} 
& 	 	x_1 \left(
\ba{c|c||c|c} 
\Lambda  & \phantom{\times} & \phantom{\times} &  \phantom{\times} \\ \hline
& & &  \\ \hline \hline 
& & &  \\ \hline
& & &   \\ 
\ena \right) & 
+ x_2 \left(
\ba{c|c||c|c} 
\ti & \ti & \ti &  M  \\ \hline
\ti & \Sigma  &  &   \\ \hline \hline 
\ti &  & - I_s &   \\ \hline
M^T & & &   \\ 
\ena \right) & \preceq & \left(\ba{c|c||c|c} 
I_p & \phantom{\ti} & \phantom{\ti} &  \phantom{\ti} \\ \hline
\phantom{\ti} & I_{r-p}  & \phantom{\ti} &  \phantom{\ti} \\  \hline  \hline  
\phantom{\ti} & \phantom{\ti} & \phantom{\ti} &  \phantom{\ti} \\ \hline 
\phantom{\ti} &  \phantom{\ti}  & \phantom{\ti} &  \phantom{\ti} \\ \hline 
\ena \right), 
\ena
\eeq
where $\Lambda$ and $\Sigma$ are diagonal, $\Lambda$ is positive definite, 
$M \neq 0, \,$  and $s \geq 0.$ 
\qed
\ecor

\section{A cookbook to generate SDPs with positive gaps}
\label{section-similar-in-anydimension} 

While two variable SDPs may  come across 
as too special, we  now show that 
they  
help us understand   positive gaps in larger SDPs: we present  three families of SDPs in which the 
same structure causes the duality gap as in the two variable case.

The SDPs in 
Examples \ref{example-E_ij} and \ref{example-E_ij_infinity} have a certain ``single sequence" structure, and they are larger versions of Example 
\ref{ex2}. To be precise, the primal optimal value is zero, 
while the dual is equivalent to a problem  like \eqref{dred} with $s=0, \,$ and therefore  has a positive optimal value.

The SDPs in Example \ref{example-double} have a richer, a certain ``double sequence" structure. 
These SDPs are  more subtle: 
we will show that the singularity degrees of {\em two} associated duals, namely of \eqref{d} and of \eqref{HD} are the largest that permit a positive duality gap.

\co{Our SDPs are inspired by the semidefinite {\em systems} in \cite{Tuncel:11}, which have a large singularity degree. 
However, as far as we can tell, our SDPs in this section are the first ones with a positive duality gap, a large singularity degree, and a 
structure that is inherited from the two variable case.  }

\subsection{Positive gap SDPs with a single sequence} 

\bex \label{example-E_ij} 
Let $n \geq 3, \, m=n-1, \,$ and let $E_{ij} \in \sym{n}$ be a matrix whose 
only nonzero entries are $1$ in positions $(i,j)$ and $(j,i).$ For brevity, let $E_i := E_{ii}.$ 

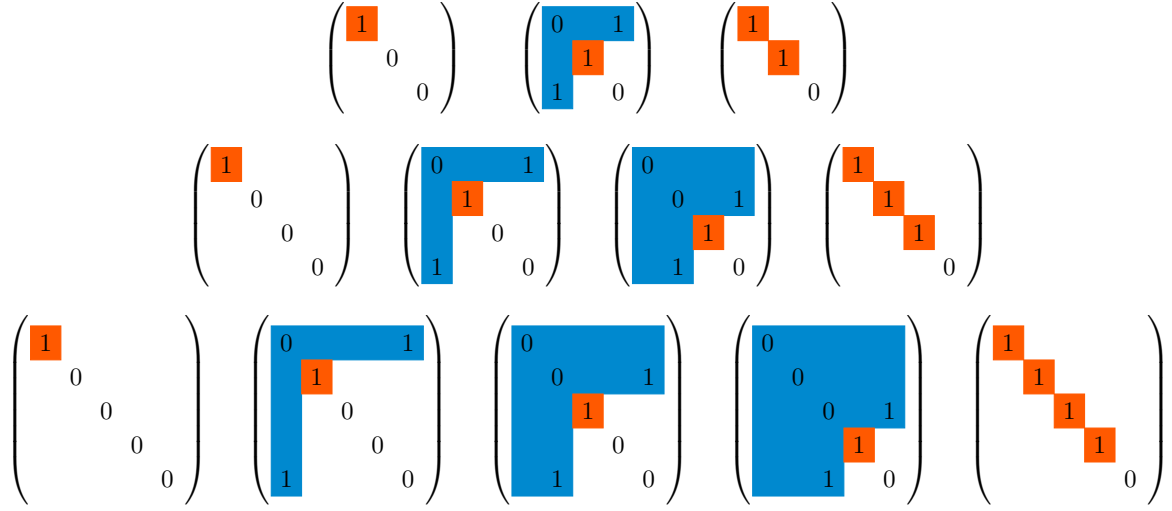
\begin{figure}[ht] 
	\begin{center} 
		\subimport{imgs/}{pic_simpleEx23v4.tex}
	\end{center} 
	\caption{The structure of the $A_i$ and of $B$  in Example \ref{example-E_ij} when $m=2, \, m=3$ and $m=4$} 
	\label{figure-sparsity-E_ij}  
\end{figure}

We  consider the SDP 
\beq \label{problem-E_ij} 
\ba{rl} 
\sup & x_{n-1}  \\
s.t.   & x_1 E_{1} +  \sum_{i=2}^{n-1} x_i (E_{i} + E_{i-1,n}) \preceq 
\bpx I_{n-1} &  \\
                   & 0 \epx. 
\ena 
\eeq 
For example, when $n=3, \, $ we recover Example \ref{ex2}.
For $n=3, 4, $ and $5$ we show the structure of the $A_i$ and of $B$ in  Figure \ref{figure-sparsity-E_ij}.
(The last matrix in each row is $B.$)

In all matrices the nonzero diagonal entries are red, and we explain the meaning of the blue submatrices shortly. 

\vspace{.1cm}  
\vspace{.1cm} 
We claim 
that there is a duality gap of $1$ between (\ref{problem-E_ij}) and its dual. 

First we compute 
the optimal value of (\ref{problem-E_ij}).
 If $x$ is feasible in 
it and $Z = (z_{ij}) \succeq 0$ is the corresponding slack, then 
$z_{nn}=0, \, $ so the last row and column of $Z$ is zero. 
Since $z_{n-2,n} = - x_{n-1}, \,$ we deduce $x_{n-1} = 0$ and 
$\val(\ref{problem-E_ij})=0.$ 

On the other hand, suppose $Y = (y_{ij})$
is feasible in the dual. By the first dual constraint 
$y_{11} = 0.$ Thus 
                \beq \label{eqn-y_ij=0} 
                \ba{rclr} 
                y_{11} = 0  & \Rightarrow & y_{1j} = 0 \, \forall j \,&  ({\rm since } \; Y \succeq 0 ) \\
                & \Rightarrow & y_{22} = 0                                & ( {\rm since}   \; A_2 \bullet Y = 0) \\ 
                & \Rightarrow & y_{2j} = 0 \, \forall j                  &   ({\rm since } \; Y \succeq 0 )  \\ 
                & \vdots & \\
                 & \Rightarrow & y_{n-2, n-2} = 0                                & ( {\rm since}   \; A_{n-2} \bullet Y = 0) \\ 
                & \Rightarrow & y_{n-2,j} = 0 \, \forall j             &  ({\rm since } \;   Y \succeq 0 ).   
                \ena
                \eeq
                We can follow this argument on Figure \ref{figure-sparsity-E_ij}.           
                If $i < m$ and $A_1 \bullet Y = \dots = A_i \bullet Y = 0,$ then the first $i$ rows 
                (and columns) of $Y \succeq 0$ are zero. Hence when we compute  $A_{i+1} \bullet Y,$ the first $i$ rows and columns $A_{i+1}$ do not  matter.
                Accordingly, we colored this portion of $A_{i+1}$   blue.

    Thus the dual is equivalent to the SDP 
    \beq \label{problem-lastreduced2}  
           \ba{rrcl}  
           \inf &   \bpx 1  & 0  \\ 0 & 0 \epx  \bullet Y'  \\ 
           s.t. & \bpx 1 & 0  \\ 0 & 0   \epx  \bullet Y'  & = & 1  \\
           & Y' & \succeq & 0,
           \ena
           \eeq
   which has optimal value $1.$ (We can think of $Y'$ as the lower right $2 \times 2$ corner of the  dual variable matrix $Y.$)
So  the duality gap is  $1,$ as wanted. 

\eex

Note the  recursive nature of the SDPs in Example \ref{example-E_ij}: if we delete the first row and column in all $A_i$ and delete $A_1, \,$ we obtain an SDP with the same structure, just with $n$ and $m$ reduced by one.


At first, these SDPs look  unnecessarily complicated, as we could just use 
$A_i = E_{i}$ for  $i=1,\dots,n-2$ and 
still have the same duality gap: the argument given in Example \ref{example-E_ij} would carry over verbatim. 
 However, this simpler  SDP is not a ``bona fide" $n-1$ variable SDP, 
since  we could simplify it even more: we could  replace $A_{1}$ with 
 $A_1 + \dots + A_{n-2}$ and drop $A_2, \dots, A_{n-2}$ 
 to obtain  a two variable SDP 
 $$
 \ba{rl}
 \sup & x_2 \\
 s.t.   & x_1 \bpx I_{n-2} & 0 \\
                           0   & 0 \epx      + x_2 \bigl( E_{n-1} + E_{n-2,n} \bigr) \preceq \bpx I_{n-1} & 0 \\ 0 & 0 \epx 
 \ena
 $$
 with the same duality gap.

Why is such a replacement impossible in \eqref{problem-E_ij}? The $(i-1,n)$ element of all $A_i$ is nonzero, so it is not hard to check  that the only psd linear combinations of the $A_i$ are nonnegative multiples of $A_1.\, $ 
We can say more: as we show in Section \ref{section-singdegree}, the matrices
$A_1, \dots, A_{n-2}$ are a minimal sequence (in a well defined sense) that zero out the first $n-2$ rows and columns of the  
dual variable matrix $Y.$ 
 
 Next comes a family of SDPs with {\em infinite} duality gap. 
\bex \label{example-E_ij_infinity} 
Let us change the last matrix in Example \ref{example-E_ij} 
to $-E_{n-1} + E_{n-2, n}.$ The resulting primal  SDP still has 
zero optimal value, but now the dual is equivalent to 
\beq \label{problem-lastreduced2}  
\ba{rrcl}  
\inf &   \bpx 1  & 0  \\ 0 & 0 \epx  \bullet Y'  \\ 
s.t. & \bpx -1 & 0  \\ 0 & 0   \epx  \bullet Y'  & = & 1,  \\
& Y' & \succeq & 0,
\ena
\eeq
      hence it is infeasible. Thus we have       
      $$
      0 \, = \, \val(\ref{p}) < \val(\ref{d}) = + \infty,
      $$
      i.e., an infinite duality gap. 
\eex 

We next discuss the connection of infinite duality gaps with weak infeasibility, another pernicious pathology of SDPs. 
We say that the dual
\eqref{d} is {\em weakly infeasible}, if the affine subspace 
$$
\{ \, Y \, | \, A_i \bullet Y = c_i \, (i=1, \dots, m) \, \} 
$$
has zero distance  to $\psd{n}$ but does not intersect it. Weakly infeasible SDPs are very challenging  for SDP solvers, which often 
mistake such instances for feasible ones. We refer to 
\cite{Waki:12, lourencco2016structural, liu2017exact} for   theoretical and computational studies and for instance libraries. 

The following proposition is folklore, but for the sake of completeness we give a short proof. 
\bprop \label{prop-weakinfeas} 
Suppose \eqref{p} is feasible, $\val \eqref{p} < + \infty, \, $ and \eqref{d} is infeasible. Then \eqref{d} is weakly infeasible. 
\eprop
\pf{} It is well known that \eqref{d} is weakly infeasible iff it is infeasible, and its alternative system 
\beq \label{eqn-alt-D} 
\ba{rcl} 
\sum_{i=1}^m {x}_i A_i & \preceq & 0 \\
\sum_{i=1}^m {x}_i c_i & = & 1
\ena  
\eeq
is also infeasible \footnote{The system \eqref{eqn-alt-D} is called an alternative system, since when it {\em is}  feasible, it is 
 a convenient certificate that \eqref{d} is infeasible:  a simple argument shows that both cannot be feasible.}.    
   So we assume that the conditions of our proposition are met, and we will show that \eqref{eqn-alt-D} is infeasible. Indeed, if  
   \eqref{eqn-alt-D}  {\em were} feasible, then adding a large multiple of a feasible solution of \eqref{eqn-alt-D} to a feasible solution of \eqref{p} would prove  $\val \eqref{p} = + \infty, \, $ which would be  a
contradiction.
\qed

Proposition \ref{prop-weakinfeas} tells us that  infinite duality gap in SDPs gives rise 
to weak infeasibility in all ``interesting" cases. Indeed, the other case of an infinite duality gap is when 
both \eqref{p} and \eqref{d} are infeasible; however, such instances are easy to produce even in linear programming. 

Proposition \ref{prop-weakinfeas} also implies  that the dual SDPs in  Example \ref{example-E_ij_infinity} are weakly infeasible. Interestingly, they are 
much simpler than the weakly infeasible SDPs  in 
\cite{Waki:12, liu2017exact}, 
while they are just as difficult, as we will show in  Section  \ref{section-computational}.


\subsection{Positive gap SDPs with a double sequence}

We now present another family of SDPs with a positive duality gap.         These may not be  {\em per se} more difficult than the ones in 
Examples \ref{example-E_ij}  and 
\ref{example-E_ij_infinity} (as we will see in Section \ref{section-computational}, those are already very hard). 
The SDPs in this section, however, have a more sophisticated ``double sequence" structure and 
we will show in Sections
\ref{section-singdegree} and \ref{section-maximal-deg} 
that the so-called singularity degree of 
{\em two} associated duals -- of (\ref{d}) and of \eqref{HD} -- are the 
maximum that permit a positive duality gap.

\bex \label{example-double}
Let $m \geq 2, \, n = 2m+1, \,$ and consider the SDP 
\beq \label{problem-double} 
\ba{rlcl}
\sup  & x_m \\ 
s.t. & x_1 \bigl( E_1 + E_{m+1} \bigr) + \sum_{i=2}^{m-1} x_i \bigl(    E_i + E_{m+i} + E_{i-1,n} + E_{m+i-1, n}   \bigr) \\
            & \hspace{3cm} \,\,\,\,\,\,\,\,\,\,\, + x_m \bigl(    E_m - E_{2 m} + E_{m-1,n} + E_{2m-1,  n}   \bigr) & \preceq & \bpx I_{m+1} &   \\   & 0 \epx.  
\ena 
\eeq 
Note that the negative sign of $E_{2m}$ in the last term is essential: if we change it to positive, then a simple calculation shows that the resulting SDP will have zero gap
with its dual.  

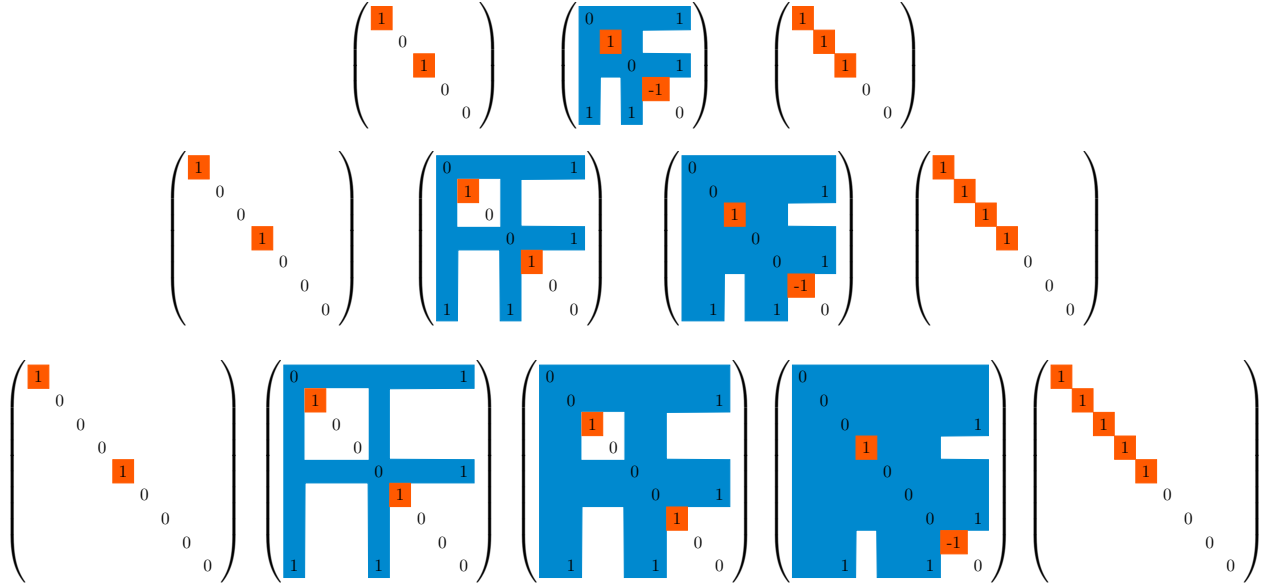
\begin{figure}[ht] 
	\begin{center} 
		\subimport{imgs/}{pic_simpleEx4v5.tex}
	\end{center} 
\caption{The structure of the $A_i$ and of $B$ in Example \ref{example-double} when $m=2, \, m=3, \,$ and $m=4$} 
\label{figure-double-sequence} 
\end{figure}

For concreteness, when $m=2$ the SDP (\ref{problem-double}) is 
\beq 
\ba{lllcl} 
\sup  & x_2 \\ \\ 
s.t. & \, x_1 \left( \begin{array}{ccc|cc} 1 & \pha{0}  & \pha{0}  & \pha{0} & \pha{0}  \\ 
	                                               \pha{0}  & 0 & \pha{0}  & \pha{0}  & \pha{0} \\ 
                                                   \pha{0} & \pha{0}  & 1 & \pha{0}  & \pha{0} \\ \hline 
                                                   \pha{0}  & \pha{0}  & \pha{0}  & 0 & \pha{0} \\
                                                   \pha{0}  & \pha{0}  & \pha{0}  & \pha{0} & 0 
                                                   \end{array} \right) & 
+ x_2 \left( \begin{array}{ccc|cc} 0 & \phao & \phao & \phao & 1 \\ 
	\phao & 1 & \phao & \phao & \phao \\ 
	 \phao & \phao & 0 & \phao & 1 \\ \hline 
      \phao & \phao & \phao & -1 & \phao \\
      1  & \phao & 1  & \phao & 0 
         \end{array} \right) & 
\preceq & \left( \begin{array}{ccc|cc} 1 & \phao & \phao & \phao & \phao \\ 
	\phao & 1 & \phao & \phao & \phao \\ 
	 	\phao & \phao & 1  & \phao & \phao \\  \hline 
        \phao & \phao & \phao  & 0 & \phao  \\ 
        \phao & \phao & \phao & \phao & 0 \end{array} \right). 
\ena
\eeq
 Figure \ref{figure-double-sequence}  depicts the structure of the $A_i$ and of $B$ in  Example 
 \ref{example-double} when $m=2,3, \, $ or $4$. The last matrix in each row is $B.$ 
 The color coding is similar to the one we used in Figure
 \ref{figure-sparsity-E_ij}, namely the nonzero diagonal elements are red, and we explain the meaning of the blue blocks  soon.  

We next argue that (\ref{p}) and (\ref{d}) satisfy 
$$0 = \val(\ref{p}) < \val(\ref{d}) = 1.$$ 

Indeed, $\val \eqref{p} = 0$ since if $x$ is feasible in (\ref{p}), then $x_m=0;$ this follows just like  in Example \ref{example-E_ij}.

Suppose next that $Y \succeq 0$ is feasible in the dual.
Then $A_1 \bullet Y = 0$ so $y_{11}=y_{m+1,m+1} = 0,$ hence 
\beq \label{eqn-y_ij=0-double} 
\ba{rclr} 
y_{11} = y_{m+1,m+1} = 0  & \Rightarrow & y_{1j} = y_{m+1,j} = 0 \, \forall j & ({\rm by} \; Y \succeq 0) \\
                 & \Rightarrow & y_{22} = y_{m+2,m+2} = 0         &  ({\rm by}  \; A_2 \bullet Y = 0)   \\   
& \Rightarrow & y_{2j} = y_{m+2,j} = 0 \, \forall j & ({\rm by} \, Y \succeq 0) \\ 
& \vdots & \\
& \Rightarrow & y_{m-1,m-1} = y_{2m-1, 2m-1} = 0 &   ({\rm by}  \;     A_{m-1} \bullet Y = 0 ) \\  
& \Rightarrow & y_{m-1,j} = y_{2m-1, j} = 0 \, \forall j &    ({\rm by} \, Y \succeq 0). \\ 
\ena
\eeq  
We can follow this argument on Figure \ref{figure-double-sequence}.  If $i < m, \, $ and $Y \succeq 0, \, $ then $A_1 \bullet Y = \dots = A_i \bullet Y = 0$ implies that the rows and columns of $Y$ indexed by $1, \dots, i$ and $m+1, \dots, m+i$ are zero. Consequently, 
these rows and columns of $A_{i+1}$ do not matter when we compute $A_{i+1} \bullet Y,$  so we colored them  blue. 

We show a feasible $Y \succeq 0$ and $A_m$  in equation 
(\ref{equation-Y-Am}) below.   (As always, the empty blocks are zero, and the $\ti$ blocks may be nonzero.)

\beq \label{equation-Y-Am}    
\underbrace{\begin{pmatrix}[c|c|c|c|c]
		\bovermat{$ m-1  $}{\phantom{0000}}	& \bovermat{$1$}{\phantom{0000}}	& \bovermat{$m-1 $}{\phantom{0000}}	& \bovermat{$1 $}{\phantom{0000}}	& \bovermat{$1 $}{\phantom{0000}}\\ \hline
		\pha{0} & \ti  & \pha{0} & \ti & \ti \\ \hline
		\pha{0} & \pha{0} & \pha{0} & \pha{0}  & \pha{0} \\ \hline
		\pha{0} & \ti & \pha{0} & \ti & \ti \\  \hline 
		\pha{0} & \ti & \pha{0} & \ti & \ti 
	\end{pmatrix}}_{Y} 
\underbrace{\begin{pmatrix}[c|c|c|c|c]
	\bovermat{$m-1  $}{\phantom{0000}} 	& \bovermat{$1$}{\phantom{0000}}	& \bovermat{$m-1$}{\phantom{0000}}	& \bovermat{$1$}{\phantom{0000}}	& \bovermat{$ 1$}{ \phao \ti \phao }\\ \hline
	 & 1  & & \pha{0} & \pha{0} \\ \hline
	 &  & &   & \ti \\ \hline
	 & \pha{0} &  & -1 & \pha{0} \\ \hline 
	\ti & \pha{0} & \ti  & \pha{0} & \pha{0} 
	\end{pmatrix}}_{A_m}, \,  
\eeq 
Thus 
(\ref{d}) is equivalent to 
\beq \label{problem-lastreduced}  
\ba{rrcl}  
\inf &   \bpx 1  & 0  \\ 0 & 0 \epx  \bullet Y'  \\ 
s.t. & \bpx 1 & 0  \\ 0 & -1  \epx  \bullet Y'  & = & 1,  \\
& Y' & \succeq & 0,
\ena
\eeq
hence it has optimal value $1, \, $ as wanted. 
\eex

\section{The singularity degree of the duals of our positive gap SDPs}
\label{section-singdegree} 

We now study our positive gap SDPs in more depth. 
We introduce faces, facial reduction, and singularity degree of SDPs,
and show that the duals associated with our SDPs, 
namely (\ref{d}) and (\ref{HD}) (defined in the Introduction), 
have singularity degree equal to  $m-1$ and $m, \,$ respectively. 

We first recall that a 
set $K$ is  a {\em cone}, if $x \in K, \, \lambda \geq 0$ implies $\lambda x \in K,$ and the dual cone of cone $K$ is
$$
K^* \, = \, \{ \, y \, | \, \la y, x \ra \geq 0 \, \forall \, x \in K \, \}. 
$$
In particular, $(\psd{n})^* = \psd{n}$ with respect to the $\bullet$ inner product.

\subsection{Facial reduction and singularity degree} 
  
\begin{Definition}

Given a closed convex cone $K,$ a convex subset $F$ of $K$ is  a {\em face of} $K,$  if 
$x, y \in K, \, \frac{1}{2}(x+y) \in F$ implies $x, y \in F.$

\end{Definition}

We are mainly interested in the faces of $\psd{n},$ \co{and the minimal cone of sets 
$H \cap \psd{n}, \,$ where $H \subseteq \psd{n}$ is an affine subspace.} 
which  have a simple and attractive description:  they are 
\beq \label{eqn-psdfaces} 
F \, \, = \, \biggl\{   T \bpx X & 0 \\ 0 & 0 \epx T^T : X \in \psd{r} \biggr\}, \, {\rm with \, dual \, cone \,} \, F^* \, = \, \biggl\{   T^{-T} \bpx X & Z \\ Z^T & Y \epx T^{-1}: X \in \psd{r} \biggr\},
\eeq 
where $0 \leq r \leq n$ and $T \in \rad{n \times n}$ is invertible (see, e.g., \cite{Pataki:00A}). 

In other words, the faces are of the form $T (\psd{r} \oplus 0)T^T$ for some $r$ and for some invertible matrix $T.$

For such a face, assuming $T=I$ we sometimes use the shorthand 
\beq \label{f-oplus}
F = \begin{pmatrix}  \oplus \! &  0  \\
0   \! &  0 
\end{pmatrix}, 
\,  F^{*} = \begin{pmatrix}  \oplus &  \ti  \\
\ti   \! &  \ti 
\end{pmatrix}, \, 
\eeq
when the size of the partition is clear from the context. 
The $\oplus$ sign denotes a positive semidefinite submatrix and the sign
$\ti$ stands for a submatrix with arbitrary elements. 

\begin{figure}[ht] 
\begin{center} 
\includegraphics[width = 10cm]{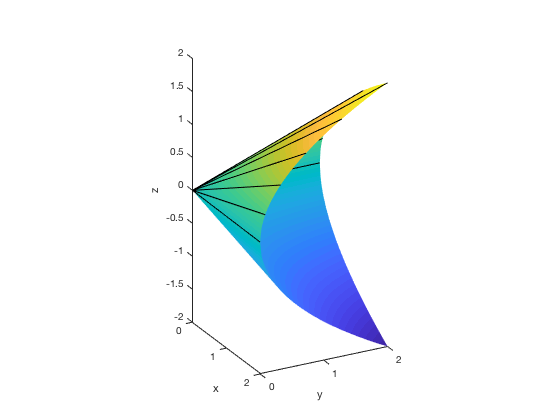}
\end{center} 
\caption{The $2 \times 2$ semidefinite cone}
\label{figure-psd} 
\end{figure}

Figure \ref{figure-psd} depicts the cone $\psd{2}$ in 3 dimensions: it plots the triplets $(x,y,z)$ such that 
\beqast
\bpx x & z \\
        z & y \epx \succeq 0.
        \eeqast
It is clear that all faces of $\psd{2}$ that 
are different from $\{ 0 \}$ and itself are extreme rays of the form 
\mbox{$ \bigl\{ \lambda uu^T : \, \lambda \geq 0 \bigr\}, $} 
where $u \in \rad{2}$ is nonzero, i.e., we can choose 
$r=1$ in (\ref{eqn-psdfaces}). 
 
\begin{Definition}
Suppose $K$ is a closed convex cone, 
$H$ is an affine subspace, and  $H \cap K \neq \emptyset.$ We define the 
{\em minimal cone of $H \cap K$} as the smallest face of $K$ that contains $H \cap K.$

We define the minimal cone of \eqref{p}, of \eqref{d} and of \eqref{HD} as   the minimal cone of their
feasible sets. In particular, the minimal cone of \eqref{p}  is the minimal cone of 
	$(\lin \{ \, A_1, \dots, A_m \} + B) \cap \psd{n}.$ 
\end{Definition}

The following easy-to-verify fact 
will help us identify the minimal cone of  SDPs: if $H \subseteq \sym{n}$ is an affine subspace which contains a psd matrix, 
\co{and 
$H \cap \psd{n} \neq \emptyset, \,$ }
then the minimal cone of $ H \cap \psd{n}$ 
is the smallest face of $\psd{n}$ that contains the maximum {\em rank} psd
matrix of $H.$ 

\bex \label{ex-mincone} 
Let $H$ be the linear subspace spanned by the matrices
$$
A_1 = \bpx 1 & \phao & \phao \\ \phao  & 0 & \phao \\ \phao & \phao  & 0 \epx, \, A_2 = \bpx \phao & \phao  & 1 \\ \phao & 1  & \phao  \\ 1 & \phao  & \phao  \epx,
$$
and $K = \psd{3}.$ 


Then 
$ H \cap K  = \{ \lambda A_1 : \lambda \geq 0 \}, \, $ hence this latter set is the minimal cone of $H \cap K.$

\eex

In this paper we are mostly interested in the minimal cone of (\ref{d}) and of (\ref{HD}) 
\footnote{In contrast, the minimal cone of (\ref{p}) is easy to identify, since it is 
	just the smallest face of $\psd{n}$ that contains the right hand side $B$ 
	(see Assumption \ref{assumption-slack}).}.

\bex \label{ex2-cont} (Example \ref{ex2} continued) 
In this example we proved that in any feasible solution of \eqref{d} 
the first row and column is zero. Thus 
\begin{equation} \label{eqn-Y} 
Y \, = \, \bpx 0 & \pha{0} & \pha{0} \\ \pha{0} & 1 & \pha{0} \\ \pha{0} & \pha{0} & 1 \epx 
\end{equation}
is a maximum rank feasible solution in (\ref{d}),  so the minimal cone of \eqref{d}   is 
$0 \oplus \psd{2}.$

\eex

Why is the minimal cone interesting? Suppose $F$ is the minimal cone of \eqref{d}. Then there is a feasible $Y$ 
in the relative interior of $F,$ otherwise the feasible set of (\ref{d}) would be contained in a smaller face of $\psd{n}.$ 
Thus replacing the primal constraint by $B - \sum_i x_i A_i \in F^*$ yields a primal-dual pair with no duality gap and primal 
attainment. (An analogous result holds for the minimal cone of (\ref{p}) and enlarging the dual feasible set).
For details, see e.g. \cite{liu2017exact}. 

How do we actually compute the minimal cone of $H \cap K$? The following basic {\em facial reduction algorithm} is designed for  this task.

\begin{algorithm}[H] 
\caption{Facial Reduction}
\label{algo:FRA}
\begin{algorithmic}
\State Let $F_0 = K, \, i=1.$
\For{${i=1,2, \dots}$}  
\State Choose $y_i \in F_{i-1}^* \cap H^\perp.$ 
\State Let $F_i = F_{i-1} \cap y_i^\perp.$
\EndFor 
\end{algorithmic}
\end{algorithm}
\begin{Definition} \label{defn-frs} 
We say that a sequence $( y_1, \dots, y_k)$ output by Algorithm 
\ref{algo:FRA} is a {\em facial reduction sequence for $K;$} 
 and we say that it is a 
{\em strict facial reduction sequence for $K$} if in addition 
 $y_i \in F_{i-1}^{*} \setminus F_{i-1}^{\perp}$ 
 for all $i \geq 1.$
 
 We denote the set of facial reduction sequences for $K$ by $\fr(K).$ 
\end{Definition}

Suppose that $(y_1, \dots, y_k)$ is constructed by Algorithm \ref{algo:FRA}. Then all $F_i$  contain
$H \cap K, \,$ and hence they also contain the 
minimal cone of $H \cap K.$ Further, 
 there is always a possible output $(y_1, \dots, y_k)$ 
such that $F_k$ is  the minimal cone, see e.g., \cite{liu2017exact}. 
We will say that such a sequence {\em defines} the minimal cone of 
$H \cap K.$ 

Clearly, Algorithm \ref{algo:FRA} can generate many possible sequences
(it can even choose several $y_i$ which are zero), but it is best  to terminate it in  a minimim number of steps. 

\begin{Definition}
	Suppose $H$ is an affine subspace with $H \cap K \neq \emptyset.$ 
The {\em singularity degree} $d(H \cap K)$ of $H \cap K$ is the smallest number of 
steps necessary for Algorithm \ref{algo:FRA} to construct the minimal cone of $H \cap K.$ 

The singularity degrees of \eqref{p}, \eqref{d}, and \eqref{HD} are defined as the singularity degree of their feasible sets.
They are denoted by $d \eqref{p}, \, d \eqref{d}, \, $ and $d \eqref{HD}, \, $ respectively.
\end{Definition}

The singularity degree of SDPs was introduced in the seminal paper \cite{Sturm:00}. 
It was used to bound the distance of a symmetric matrix  from $H \cap \psd{n}, \,$ given 
the distances from $H$ and from $\psd{n}.$ More recently it was used 
in \cite{drusvyatskiy2017note} to bound the rate of convergence of the alternating 
projection algorithm to such a set. 

For later reference, we state a basic bound on the singularity degree (for details, see, e.g., \cite[Theorem 1]{liu2017exact}):  
\begin{equation} \label{eqn-sing-Hperp} 
d(H \cap K) \leq \dim H^\perp. 
\end{equation}

In the following examples involving SDPs 
we denote the members of facial reduction sequences by capital letters (since they are matrices).

\bex (Example \ref{ex-mincone} continued)    
In this example the sequence $(Y_1, Y_2)$ below defines the minimal cone, with 
corresponding faces shown. Note that $F_2$ is the minimal cone.

\beqast
Y_1 & =  & \bpx 0 & 0 & 0 \\ 0 & 0 & 0 \\ 0 & 0 & 1 \epx, \, F_1 \, = \, \psd{3} \cap Y_1^\perp \, = \, \bpx 
\parbox{1cm}{$\,\, \bigoplus$}   & \hspace{-.4cm} \parbox{0.4cm} {$\ba{c} 0 \, \\ 0 \, \ena$} \hspace{-5.9cm} \\ 
\parbox{1cm}{$0 \,\,\,\,\,\, 0$}  &   \hspace{-.3cm} 0 \hspace{0.4cm} 
\epx,  \\
Y_2 & = & \bpx 0 & 0 & -1 \\ 0 & 2 & 0 \\ -1 & 0 & 0 \epx, F_2 \, = \, F_1 \cap Y_2^\perp \, = \, \bpx \oplus & 0 & 0 \\
0 & 0 & 0 \\
0 & 0 & 0 
\epx.
\eeqast 

Since $Y_1$ is the maximum rank psd matrix in 
$H^\perp, \,$ it is the best choice to start such a facial reduction sequence, hence 
$d(H \cap \psd{3}) = 2.$ 

\eex

\bex (Example \ref{ex2} and \ref{ex2-cont} continued) In this example the minimal cone of \eqref{d} is $0 \oplus \psd{2}$ and 
$A_1 \bullet Y = 0$ implies $Y \in 0 \oplus \psd{2}.$
Since $A_1 \succeq 0, \,$ 
the one element facial reduction sequence 
$(A_1)$ defines the minimal cone of \eqref{d}. 
\eex

The reader may wonder, why we connect positive gaps to the singularity degree of (\ref{d}) and of $(HD), \,$ and not to the
singularity degree of 
$(P).$ We could  do the latter, by exchanging  the roles of the 
primal and dual. However, we think that our  treatment is more intuitive, as we next explain.

The dual feasible set is 
$H \cap \psd{n}, \,$ where 
\beq \label{eqn-H-for-D} 
\begin{array}{rcl}
H & = & \{ \, Y \in \sym{n} \, | \, A_i \bullet Y = c_i \, \forall i \, \} \\
    & = &  \{ \, Y \in \sym{n} \,  | \, A_i \bullet Y = 0 \, \forall i \, \} + Y_0,   
\end{array}
\eeq 
where $Y_0 \in H$ is arbitrary. 
Thus, to define the minimal cone of (\ref{d}) we use a facial reduction sequence whose members are in $H^\perp \subseteq \lin \{ A_1, \dots, A_m \ \}.$
As we will show, in our instances actually the $A_i$ themselves form a facial reduction sequence 
that defines the minimal cone of \eqref{d}, and this makes the essential structure of the minimal cone apparent.  
An analogous statement holds for the minimal cone of $(HD).$ 

\subsection{The singularity degree of the single sequence SDPs in Example \ref{example-E_ij}} 

We now analyze the singularity degree of the duals of the SDPs given in Example \ref{example-E_ij}.
\bth 
\label{thm-singdegree-single} 
Let (\ref{d}) be the dual of the SDP (\ref{problem-E_ij}).
Then 
$$
d(D) = m-1.
$$
\enth 
\pf{} Recall that $m=n-1$ in this SDP. We first claim that the minimal cone of \eqref{d} is $0 \oplus \psd{2}.$ Indeed, by 
the argument in (\ref{eqn-y_ij=0})  the minimal cone is contained in $0 \oplus \psd{2}.$
Since $0 \oplus I_2$ is feasible in \eqref{d}, the minimal cone is exactly $0 \oplus \psd{2}.$

An analogous argument (by plugging $Y \succeq 0$ into the equations $A_1 \bullet Y = 0, \, \dots, \, A_{i} \bullet Y = 0$, like in (\ref{eqn-y_ij=0})) proves 
\beq \label{eqn-psd-Ai} 
\psd{n} \cap A_1^\perp \cap \dots \cap A_i^\perp \, = \, 0_i \oplus \psd{n-i} \; {\rm for} \; 1 \leq i \leq m-1. 
\eeq 
Let $i \in \{1, \dots, m-1 \}.$ Then clearly $A_{i+1} \in ( 0_i \oplus \psd{n-i})^*, \, $ so by  \eqref{eqn-psd-Ai} we deduce 
$$
A_{i+1} \in ( \psd{n} \cap A_1^\perp \cap \dots \cap A_i^\perp)^*, 
$$
hence $(A_1, \dots, A_{m-1})$ is a facial reduction sequence, which defines the minimal cone of 
$(D). \,$ Thus $d(D) \leq m-1.$ 

To complete the analysis we  show that any strict facial reduction sequence in 	
$\lin \, \{ A_1, \dots, A_m \}$ reduces $\psd{n}$ by at most as much as the $A_i$ themselves. 
This is done in Claim \ref{Claim-D}, whose proof   is in Appendix \ref{appendix-proof-Yi-Ai}.
 
\begin{Claim} \label{Claim-D} 
Suppose $i \in \{1, \dots, m-1 \}$ and $(Y_1, \dots, Y_i)$ is a strict facial reduction sequence, whose members are all in 
$\lin \{ A_1, \dots, A_m \}.$ Then 
\beq \label{eqn-Yi-Ai} 
\psd{n} \cap Y_1^\perp \cap \dots \cap Y_i^\perp \, = \, 0_i  \oplus \psd{n- i}. 
\eeq 
\end{Claim} 

Having Claim \ref{Claim-D} at hand, let $(Y_1, \dots, Y_i)$ be  a strict 
facial reduction sequence that defines the minimal cone of \eqref{d}. 
Since the minimal cone is $0_{n-2} \oplus \psd{2},$ by Claim 
\ref{Claim-D} we deduce  $i = n-2 = m-1.$ Hence 
$d(\ref{d}) = m-1,$ and the desired result follows. 
\qed

\subsection{The singularity degree of the double sequence SDPs in Example \ref{example-double}} 

We now turn to studying the singularity degrees of the duals of the SDPs in Example 
\ref{example-double}. 
The main rationale for creating these SDPs in the first place is that the singularity degrees of {\em two} associated duals achieve an upper bound. 
\bth 
\label{thm-singdegree-double} 
Let (\ref{d}) be the dual of (\ref{problem-double}) and \eqref{HD}  its homogeneous dual.
 Then 
$$
d \eqref{d}  = m-1 \,\, {\rm and} \,\, d \eqref{HD}  = m.
$$
\enth 
\pf{sketch} 
Recall that $n=2m+1$ in this example. 
The proof of $d \eqref{d}  = m-1$ is almost verbatim the same as the proof of
the same statement in Theorem \ref{thm-singdegree-single}: the key is 
that for $i \in \{1, \dots, m-1 \}$ 
\beq \nonumber 
\ba{rcl}
\psd{n} \cap A_1^\perp \cap \dots \cap A_i^\perp & = & \{ \; Y \succeq 0 \; | \; {\rm \; rows \; and \; columns \; of \; } Y {\rm \; indexed \; by \; } \\
&    & \;\;\;\;\;\;\;\;\;\;\;\;\;\;\;\;\; 1, \dots, i {\rm \; and \;} m+1, \dots, m+i \; {\rm are \; zero \;} \} 
\ena
\eeq 
(cf. \eqref{eqn-y_ij=0-double}). We leave the details to the reader.

We next outline the proof of $d \eqref{HD} = m.$ The argument goes like this:
when we check whether $Y \succeq 0$ is feasible in 
\eqref{HD}, it suffices (and  is convenient)  to plug $Y$ into the equations 
$$
B \bullet Y = 0, \, A_2 \bullet Y = 0, \dots, A_m \bullet Y = 0
$$
in this order. Indeed, $B \bullet Y = 0$ implies that the first $m+1$ rows and columns of $Y$ are zero, and this implies $A_1 \bullet Y = 0,$ so we do not have to check this last equation separately.  
Then $A_2 \bullet Y = 0$ implies that the $(m+2)$nd  row and column of $Y$ is zero, and so on.

Thus, 
\begin{equation} \label{eqn-mincone-HD} 
\psd{n} \cap B^\perp \cap A_2^\perp \cap \dots \cap A_i^\perp \, = \, 0_{m+i} \oplus \psd{m-i+1} \,\, {\rm for \,} i=2, \dots, m.
\end{equation} 
Equation  \eqref{eqn-mincone-HD} with $i=m$  implies that the minimal cone of \eqref{HD} is contained in $0 \oplus \psd{1}.$   
Since $0 \oplus I_1$ is feasible in \eqref{HD}, the minimal cone is exactly $0 \oplus \psd{1}.$ 

Equation  \eqref{eqn-mincone-HD}  also implies that $(B, A_2, \dots, A_{m-1}, - A_m)$ is a facial reduction sequence that defines the minimal cone of \eqref{HD}. 
Consequently, $d \eqref{HD} \leq m.$ 

We can prove $d \eqref{HD} = m$ by using Claim \ref{Claim-HD} below, 
whose proof is very similar to the 
proof of Claim \ref{Claim-D},  so it is omitted. 

\begin{Claim} \label{Claim-HD} 
	Suppose $i \in \{1, \dots, m \}$ and $(Y_1, \dots, Y_i) $ is a strict facial reduction sequence, 
	whose members are all in 
	$\lin \, \{ \,  B, A_2, \dots, A_m \, \}.$ 
	Then  
	\beq \nonumber \label{eqn-YiAiB} 
	\psd{n} \cap Y_1^\perp \cap \dots \cap Y_i^\perp \, = \, 0_{m+i} \oplus \psd{m-i+1}.
	\eeq
\end{Claim}  

Given Claim \ref{Claim-HD} the proof of Theorem \ref{thm-singdegree-double} is complete.
\qed

\section{Maximal singularity degree implies zero duality gap} 
\label{section-maximal-deg} 

We now show that the SDPs of Section \ref{section-similar-in-anydimension} are, in a well defined sense, the best possible:
we prove that 
$$
d \eqref{d} \leq m, \, {\rm and} \, d \eqref{HD} \leq m+1 \,$$ always hold, 
and when these upper bounds are attained, there is no gap.
For the reader's sake we present our results in two subsections.

\subsection{Maximal singularity degree of \eqref{d}  implies zero duality gap} 

The main result of this subsection is fairly straightforward. 
\bprop
The inequality 
\beqast
d \eqref{d} & \leq & m.
\eeqast
always holds, and when $d \eqref{d}  = m \,$ there is no duality gap.
\eprop
\pf{} 
Let $H$ and $Y_0$ be as in \eqref{eqn-H-for-D}. We see that 
\begin{equation} \label{eqn-Hperp} 
H^\perp \, = \, \lin \{ \, A_1, \dots, A_m \, \} \cap Y_0^\perp.
\end{equation}
Since the feasible set of \eqref{d} is $H \cap \psd{n}, \, $ we deduce 
$$
d \eqref{d}  \leq \dim H^\perp \, \leq \, m, 
$$
where the first inequality comes from \eqref{eqn-sing-Hperp} and the second from \eqref{eqn-Hperp}.

Suppose $d \eqref{d}=m.$ Then $\dim H^\perp = m, \,$ hence $c_i = A_i \bullet Y_0 = 0 \,$ for all $i. \,$ So the objective 
function in \eqref{p} is $c = 0, \, $ and this proves that there is no gap.
\qed  

\subsection{Maximal singularity degree of \eqref{HD}  implies zero duality gap} 

We first observe that \eqref{eqn-sing-Hperp} with $H = \{ \, Y \, : \, B \bullet Y = A_i \bullet Y = 0 \, \forall i \, \}$ implies 
\begin{equation}
d \eqref{HD} \leq m+1. 
\end{equation}

The main result of this subsection is 
\bth \label{theorem-d-equals-mplus1}
Suppose  $d \eqref{HD}  = m+1. \,$ Then the only feasible solution of \eqref{p} is $x=0, \,$ and 
\eqref{d} is strictly feasible. Hence 
$$
\val(\ref{p}) = \val(\ref{d}) = 0. 
$$
\enth
\qed

To prove Theorem \ref{theorem-d-equals-mplus1},   we first we define certain structured facial reduction sequences for $\psd{n}.$ These sequences were originally introduced in \cite{liu2017exact}. 

\begin{Definition} \label{definition-regfr} 
	We say that $(M_1, \dots, M_k)$ is a {\em regularized facial reduction sequence for $\psd{n}$} 
	if $M_i$ is of the form 
	$$
	M_i  = 
	\bordermatrix{
		& \overbrace{\qquad \qquad \qquad}^{\textstyle r_{1}+\ldots+r_{i-1}} & \overbrace{\qquad}^{\textstyle r_{i}} & \overbrace{\qquad \qquad \qquad\quad}^{\textstyle n-r_{1}-\ldots-r_{i}} \cr\\
		& \times  &  \times  &  \times \cr
		& \times  &  I   &  0 \cr
		& \times  &  0  &  0 \cr} \in \sym{n}
	$$
	for  $i= 1,\dots, k, \,$ where the $r_i$ are nonnegative integers, and the $\times$ symbols 
	correspond to blocks with arbitrary elements.
	
	\co{We denote the set of such sequences by $\regfr(\psd{n}).$ 
	Sometimes we  refer to such a sequence 
	by the length $r_i$ blocks
	\beqast
	\I_1 & := & \{ 1, \dots, r_1 \}, \\
	\I_2 & := & \{ r_1 + 1, \dots, r_1 + r_2 \}, \\
	& \vdots & \\
	\I_k & := & \{  \sum_{i=1}^{k-1} r_i + 1,  \dots, \sum_{i=1}^{k} r_i  \}.  
	\eeqast
} 
\end{Definition}
For instance,  $(A_1, \dots, A_m)$   in the single sequence SDPs in 
Example \ref{example-E_ij} is a regularized facial reduction sequence. We refer to 
Figure \ref{figure-sparsity-E_ij} in which the identity blocks in the $A_i$ are red, and the blocks with arbitrary elements are blue. 
Also, $(B, A_2, \dots, A_{m-1}, - A_m)$ is a regularized facial reduction sequence in the double sequence SDPs in 
Example \ref{example-double}.

\ble \label{lemma-double-ref} 
Suppose $d \eqref{HD} = m +1.$ Then \eqref{p} has a reformulation 
	\beq \label{Pprime} \tag{$P'$}                                                                          
	\begin{array}{rcl}        
		\sup   \sum_{i=1}^m c_i' x_i     & &        \\                                                                             
		\sum_{i=1}^m x_i A_i'  & \preceq & B \, = \, \bpx I_r & 0 \\ 0 & 0 \epx 
	\end{array}                                                                                              
	\eeq
	in which 
$(B, A_1', \dots, A_{m}')$ is a regularized facial reduction sequence. 
		 Furthermore,   \eqref{Pprime} can be constructed using only operations 
		 \eqref{exch}, \eqref{trans} and \eqref{rotate} in Definition \ref{definition-reform}. 
\ele 

The proof of Lemma 
\ref{lemma-double-ref}  is a bit technical, so we give  in Appendix \ref{subsection-proof-thm-doubleref}. However, we next illustrate it.
\co{First, using it we can easily visualize the  maximum rank slack both in \eqref{p} and \eqref{Pprime} (as it is just the right hand side).
Second, we can visualize the maximum rank feasible solution in \eqref{HD'}, the homogeneous dual of \eqref{Pprime}. 
}
\co{As to the latter point, suppose $Y$ is feasible in the homogeneous dual of \eqref{Pprime}, and
	also suppose the block sizes in the sequence $(B, A_1', \dots, A_{d(HD)-1}')$   are $r_0, r_1, r_2, \dots.$ 
	Since $B \bullet Y = 0$ the first $r_0$ rows and columns of $Y$ are zero. 
	Since $A_1 \bullet Y = 0$ the next $r_1$ rows and columns of $Y$ are zero, and so on.
}
\co{
As to the latter point, 
suppose $Y$ is feasible in the homogeneous dual of \eqref{Pprime}, and
also suppose the block sizes in the sequence $(B, A_1', \dots, A_{d-1}')$   are $r_0, r_1, \dots r_d.$ Let $r = \sum_{i=0}^d r_i.$ 
Since $B \bullet Y = 0, \, $ the first $r_0$ rows and columns of $Y$ are zero. 
Since $A_1 \bullet Y = 0$ the next $r_1$ rows and columns of $Y$ are zero, and so on. This argument, together with some 
$Y \in 0  \oplus \pd{n-r}$ which is feasible in \eqref{HD}, proves that the minimal cone of 
\eqref{HD} is  $0 \oplus \psd{n-r}.$ 
}

\bex \label{example-doublereformulation} 
Consider the SDP 
\beq \label{eqn-doublereformulated} 
\ba{llcl}
\sup & 13 x_1  - 3 x_2 \\
s.t. & x_1 \bpx 
0 & & 2  & 2    \\
& 1 & & \\
2   &  & 0  & \\
2       &  & & 0  
\epx
+ x_2 \bpx
0  &    &    &   \\
& 0 &    & 2  \\
&    & 1 &    \\
& 2  &  & 0  
\epx 
\preceq
\bpx 
1 &    &   & \\
&  0 &   & \\
&     & 0 & \\
&     &    & 0 
\epx,
\ena 
\eeq
and let $A_1$ and $A_2$ denote the matrices on the left hand side, and $B$ the right hand side. Then 
\bit
\item $B$ is the maximum rank slack and \eqref{eqn-doublereformulated} is in the form of \eqref{Pprime}.  

\item $(B, A_1, A_2)$ is a regularized facial reduction sequence, which 
defines the
minimal cone of $(HD), \,$ which is $0 \oplus \psd{1}, \,$ i.e., the set of nonnegative multiples of 
$$
Y = \bpx 0 & 0 & 0 & 0 \\
0 & 0 & 0 & 0 \\
0 & 0 & 0 & 0 \\
0 & 0 & 0 & 1
\epx. 
$$
\eit 
Thus $d(\mbox{\ref{HD}}) \leq 3$ and we claim 
that actually $$d(\mbox{\ref{HD}}) = 3 \, {\rm holds.}$$ 

Indeed, let 
$L := \lin \, \{ B, A_1, \dots, A_m \},$ then the facial reduction sequences that define the minimal cone of (\ref{HD}) are in $L.$ 

Clearly, $B$ is the only nonzero psd matrix in $L.$ Further,   
$A_1$ is the only matrix in \mbox{$L \cap (\psd{3} \cap B^\perp)^*$}  
whose lower right 
$3 \times 3$ block is nonzero, thus 
$(B, A_1)$ is the only strict length two facial reduction sequence in $L.$ 
By similar logic  $(B, A_1, A_2)$ is the only  strict length three facial reduction sequence in $L.$ 

Thus $d \eqref{HD} =3,$ as desired.

\eex 

We need some more notation: for $Y \in \sym{n}$ and $\I, \J  \subseteq \{1, \dots, n \}$ we define 
$
Y(\I, \J)
$
as the submatrix of $Y$ with rows in $\I$ and columns in $\J,$ 
and let 
\beqast
Y(\I) & := & Y(\I, \I). 
\eeqast

\pf{of Theorem \ref{theorem-d-equals-mplus1}:} 
Assume $d \eqref{HD} =m+1.$ Clearly, we can assume that \eqref{p} was reformulated into the form of \eqref{Pprime}: as Lemma \ref{lemma-double-ref} shows, this can be done using only operations \eqref{exch}, \eqref{trans} and \eqref{rotate} in Definition \ref{definition-reform}, and under these operations the statements of Theorem \ref{theorem-d-equals-mplus1} are invariant.  

 So we assume that \eqref{p} is the same as \eqref{Pprime}, and we denote the constraint matrices on the left by $A_i$ for $i=1, \dots, m.$ For brevity, let $A_0 := B.$  

Assume that in \eqref{p} the regularized facial reduction sequence $(A_0, A_1, \dots, A_{m})$ has block sizes 
$r_0, r_1, \dots, r_m.$ Define the index sets 
	\beqast
	\I_0 & := & \{ 1, \dots, r_0 \}, \\
	\I_1 & := & \{ r_0 + 1, \dots, r_0 + r_1 \}, \\
	& \vdots & \\
	\I_m & := & \{  \sum_{i=1}^{m-1} r_i + 1,  \dots, \sum_{i=1}^{m} r_i  \}. \\
	\eeqast
Further, for $i \leq j$ we let 
\beqast
\I_{i:j} & := & \I_i \cup \dots \cup \I_j, 
\eeqast
and $\I_{m+1} := \{1, \dots, n \} \setminus \I_{1:m}.$ Finally, 
we write $\I_{i:}$ for $\I_{i:(m+1)}$ for all $i$ 
(we can do this without confusion, since $m+1$ is the largest index). 

We first prove 
\beq \label{eqn-nonzeroblock} 
A_i( \I_{i-1}, \I_{(i+1):})    \neq 0 \,\,  {\rm for} \,\, i=1, \dots, m.
\eeq 
Let $i \in \{1, \dots, m \}$ and let us picture $A_{i-1}$ and $A_i$ in equation \eqref{bla}: as always, 
the empty blocks are zero, and the $\ti$ blocks are arbitrary. The blocks marked by $\otimes$ are 
$A_i( \I_{i-1}, \I_{(i+1):})$ and its symmetric counterpart.   

\beq \label{bla}  
A_{i-1} \, = \, \begin{pmatrix}[c|c|c|c]
\bovermat{$\I_{0:(i-2)}$}{\mbox{$\,\,\,\, \times \,\,\,\,$}} 	& \bovermat{$\I_{i-1}$}{\mbox{$\,\, \times\,\,$}}	& \bovermat{$ \I_i$}{\mbox{$\,\,\times\,\,$}}	& \bovermat{$  \I_{(i+1):}$}{\mbox{$\,\,\,\,\times\,\,\,\,\,$}}	\\ \hline 
\ti & I &  &    \\ \hline
\ti &  & &    \\ \hline
\ti & \pha{0} &  &  
\end{pmatrix}, \,  
A_{i}  \, = \, \begin{pmatrix}[c|c|c|c]
\bovermat{$\I_{0:(i-2)}$}{\mbox{$\,\,\,\, \times \,\,\,\,$}} 	& \bovermat{$\I_{i-1}$}{\mbox{$\,\, \times\,\,$}}	& \bovermat{$ \I_i$}{\mbox{$\,\,\times\,\,$}}	& \bovermat{$  \I_{(i+1):}$}{\mbox{$\,\,\,\,\times\,\,\,\,\,$}}	\\ \hline 
\ti & \ti & \ti & \otimes  \\ \hline
\ti & \ti & I &    \\ \hline
\ti & \otimes &  &  
\end{pmatrix}. 
\eeq 
Now suppose the $\otimes$ blocks are  zero and let 
 $A_{i-1}' := \lambda A_{i-1} + A_i$ for some large $\lambda >0.$ Then by the Schur complement condition for positive definiteness we find 
$$
A_{i-1}'(\I_{(i-1):i}) \succ 0,
$$
hence $(A_1, \dots, A_{i-2}, A_{i-1}', A_{i+1}, \dots, A_m)$ is a shorter facial reduction sequence, 
which also defines the minimal cone of \eqref{HD}. Thus $d \eqref{HD} \leq m, \, $ which is a contradiction.
We thus proved \eqref{eqn-nonzeroblock}. 

We illustrate statement \eqref{eqn-nonzeroblock} in  equation (\ref{eqn-otimes}),  when $m=2$ and $d \eqref{HD} =3:$ the $\otimes$ blocks in $A_1$ and $A_2$ 
were just proven to be nonzero. 

\beq \label{eqn-otimes}    
A_0 = \begin{tabular}{|m{14pt}|m{14pt}|m{14pt}|m{14pt}|}
\multicolumn{1}{c}{\tikzmark{a}{\makebox[14pt][c]{\tiny{\color{white}0}}}}
&\multicolumn{1}{c}{\tikzmark{b}{\makebox[14pt][c]{\tiny{\color{white}0}}}}
&\multicolumn{1}{c}{\tikzmark{c}{\makebox[14pt][c]{\tiny{\color{white}0}}}}
&\multicolumn{1}{c}{\tikzmark{d}{\makebox[14pt][c]{\tiny{\color{white}0}}}}\\
\hline
\multicolumn{1}{|c|}{$I$}	&	&	&\\
\hline
&	&	&\\
\hline
&	&	&\\
\hline
&	&	&\\
\hline
\end{tabular}
\link{a}{a}{$\I_0$}
\link{b}{b}{$\I_1$}
\link{c}{c}{$\I_2$}
\link{d}{d}{$\I_3$} \, , \,\,\, A_1 = \begin{tabular}{|m{14pt}|m{14pt}|m{14pt}|m{14pt}|}
\multicolumn{1}{c}{\tikzmark{a}{\makebox[14pt][c]{\tiny{\color{white}0}}}}
&\multicolumn{1}{c}{\tikzmark{b}{\makebox[14pt][c]{\tiny{\color{white}0}}}}
&\multicolumn{1}{c}{\tikzmark{c}{\makebox[17pt][c]{\tiny{\color{white}0}}}}
&\multicolumn{1}{c}{\tikzmark{d}{\makebox[14pt][c]{\tiny{\color{white}0}}}}\\
\hline
\multicolumn{1}{|c|}{$\times$}  &	\multicolumn{1}{c|}{$\times$}  & \multicolumn{2}{c|}{$\otimes$}\\
\hline
\multicolumn{1}{|c|}{$\times$}  & \multicolumn{1}{c|}{$I$}	&	&\\
\hline
\multirow{2}{*}{$\otimes$}	&	&	&\\
\cline{2-4}
&	&	&\\
\hline
\end{tabular}
\link{a}{a}{$\I_0$}
\link{b}{b}{$\I_1$}
\link{c}{c}{$\I_2$}
\link{d}{d}{$\I_3$}\,,  \,\,\, A_2 =  \begin{tabular}{|m{14pt}|m{14pt}|m{14pt}|m{14pt}|}
\multicolumn{1}{c}{\tikzmark{a}{\makebox[14pt][c]{\tiny{\color{white}0}}}}
&\multicolumn{1}{c}{\tikzmark{b}{\makebox[14pt][c]{\tiny{\color{white}0}}}}
&\multicolumn{1}{c}{\tikzmark{c}{\makebox[17pt][c]{\tiny{\color{white}0}}}}
&\multicolumn{1}{c}{\tikzmark{d}{\makebox[14pt][c]{\tiny{\color{white}0}}}}\\
\hline
\multicolumn{1}{|c|}{$\times$}	& 	\multicolumn{1}{c|}{$\times$}	&		\multicolumn{1}{c|}{$\times$} & 	\multicolumn{1}{c|}{$\times$} \\
\hline
\multicolumn{1}{|c|}{$\times$} &		\multicolumn{1}{c|}{$\times$} &		\multicolumn{1}{c|}{$\times$} & 	\multicolumn{1}{c|}{$\otimes$} \\
\hline
\multicolumn{1}{|c|}{$\times$}  &		\multicolumn{1}{c|}{$\times$}  &		\multicolumn{1}{c|}{$I$}  &\\
\hline
\multicolumn{1}{|c|}{$\times$} 	&		\multicolumn{1}{c|}{$\otimes$} &	&\\
\hline
\end{tabular}
\link{a}{a}{$\I_0$}
\link{b}{b}{$\I_1$}
\link{c}{c}{$\I_2$}
\link{d}{d}{$\I_3$}
  \eeq 
  
We next prove that the only feasible solution of \eqref{p} is $x = 0.$ For that, let $x$ be feasible in \eqref{p} and define 
$$Z = A_0 - \sum_{i=1}^m x_i A_i.$$ 
Since $A_i(\I_{m+1})=0$ for all $i,\, $ we deduce $Z(\I_{m+1})=0.$ Since $Z \succeq 0, \, $ we deduce that the columns of 
$Z$ corresponding to $\I_{m+1}$ are zero. Hence
$$
0 \, = \, Z(\I_{m-1}, \I_{m+1}) \, = \, x_m \underbrace{A_m(\I_{m-1}, \I_{m+1})}_{\neq 0}, 
$$
so  $x_m=0,$ which in turn implies $Z(\I_{m:(m+1)}) =0.$ By a similar reasoning, 
$$
0 \, = \, Z(\I_{m-2}, \I_{m:(m+1)}) \, = \, x_{m-1}  \underbrace{A_{m-1}(\I_{m-2}, \I_{m:(m+1)})}_{\neq 0}, 
$$
hence $x_{m-1} = 0,$ and so on. Thus $x=0$ follows, as wanted.

We finally prove that \eqref{d}  is strictly feasible. By condition \eqref{eqn-nonzeroblock} the $A_i$ are linearly independent, 
so there exists $\bar{Y} \in \sym{n}$ such that $A_i \bullet \bar{Y} = c_i \, $ for all $i.$ 

By an argument analogous to the previuous one, we  see that the only psd linear combination of the $A_i$ is the zero matrix. 
Thus the Gordan-Stiemke theorem \eqref{eqn-HKlemma} with $$H := \lin \, \{ \, A_1, \dots, A_m \, \}$$  
implies that there is $\hat{Y} \succ 0$ such that 
$A_i \bullet \hat{Y} = 0$ for $i=1, \dots, m.$ 
Clearly, if $\lambda > 0$ is large enough, then $Y := \bar{Y} + \lambda \hat{Y}$ is strictly feasible in \eqref{d}. 

The proof is now complete.

\qed

\bex (Example \ref{example-doublereformulation} continued) 
Recall that in this SDP we have $d \eqref{HD} =3.$ 
We chose the blocks specified in equation 
(\ref{eqn-nonzeroblock}) to have entries all equal to $2.$ Note that $\I_0 = \{ 1 \}, \, \I_1 = \{2\}, \, \I_2 = \{3\}.$ 

A possible $Y$ that is strictly feasible in \eqref{d} is 
$$
Y \, = \, (y_{ij})_{i,j=1}^4 \, = \, \bpx 1 & 0 & 2  & 1 \\
0 & 1 & 0 &  -2 \\
2 & 0 & 5 & 0 \\
1 & -2 & 0 & 25 
\epx. 
$$
\eex 

It is worth  to elaborate on Theorem \ref{theorem-d-equals-mplus1} some more. 
By the Gordan-Stiemke theorem \eqref{eqn-HKlemma} the dual \eqref{d} is strictly feasible iff its singularity degree is $0.$ 
Thus Theorem \ref{theorem-d-equals-mplus1} establishes the following surprising implication:
\begin{equation}
d \eqref{HD} = m \, \Rightarrow \, d \eqref{d} = 0.
\end{equation}

\section{A computational study}
\label{section-computational}

This section presents a computational study of SDPs with positive duality gaps.

We first remark that pathological SDPs are extremely difficult to solve by interior point methods. However, 
some recent implementations of facial reduction \cite{permenter2014partial, zhu2019sieve} 
work on some pathological SDPs. We refer to  \cite{liu2019new} for 
an implementation of the Douglas-Rachford splitting algorithm to solve the weakly infeasible SDPs
from \cite{liu2019new} and to \cite{hauenstein2018numerical} 
for a homotopy method to tackle these same 
SDPs. Furthermore, the exact SDP solver SPECTRA \cite{henrion2017spectra} can solve small SDPs in exact arithmetic. 
We hope that the detailed study we present here will inspire further research. 

We generated a library of challenging SDPs based on the single sequence SDPs in Example \ref{example-E_ij}  and Example \ref{example-E_ij_infinity}. 

First we created SDPs described in Example \ref{example-E_ij}  with $m=2, \dots, 11.$ We multiplied the primal objective by $10, \, $  meaning we chose 
$$
c = 10 e_{m}.
$$
(Recall $m=n-1.$) Thus in these instances the primal optimal value is still zero, but the dual optimal value is $10.$ 

Second, we constructed single sequence SDPs as given in Example 
\ref{example-E_ij_infinity}, with $m=2, \dots, 11.$ 
 For consistency,  we multiplied the primal objective by $10,  \, $  to make it 
$10 e_{m}.$ 
Then the primal optimal value is still zero, and 
the dual is still infeasible, so the duality gap is still infinity. (Recall from Proposition \ref{prop-weakinfeas} that the dual is weakly infeasible.) 

We say that the SDPs thus created from Examples \ref{example-E_ij} and \ref{example-E_ij_infinity} 
are {\em clean}, meaning the duality gap can be verified by simple 
inspection. 

To construct  SDPs in which the duality gap is less obvious, we added an optional 
\begin{itemize}
\item[] {\bf Messing step:} Let $T$ be  an invertible matrix
with integer entries, and replace all 
$A_i$ by $T^T A_i T$ and $B$ by $T^T B T.$ 
\end{itemize}	    

Thus we have  four categories of SDPs, stored under the names 
\begin{itemize}
\item ${\mathit gap}\_{\mathit single}\_{\mathit finite}\_{\mathit clean}\_m,  \,$ 
\item ${\mathit gap}\_{\mathit single}\_{\mathit finite}\_{\mathit messy}\_m,  \,$  
\item ${\mathit gap}\_{\mathit single}\_{\mathit inf}\_{\mathit clean}\_m,$   and 
\item ${\mathit gap}\_{\mathit single}\_{\mathit inf}\_{\mathit messy}\_m,$   
\end{itemize}
where $m=2, \dots, 11$ in each category.
We tested two SDP solvers: the Mosek commercial solver, and the SDPA-GMP high precision SDP solver \cite{fujisawa2008sdpaGMP}. 

Table \ref{table1} reports the number of correctly solved instances in each category.

The solvers are not designed to  detect a finite 
positive duality gap. Thus, to be fair,  we report that  an SDP with a finite duality gap was correctly solved when 
Mosek does {\em not} report ``OPTIMAL" or ``NEAROPTIMAL" status, 
and   SDPA-GMP does {\em not} report 	  ``pdOPT" status \footnote{Of course, with such reporting we could declare a poor SDP solver to be good, just because it rarely reaches optimality. However, Mosek and SDPA-GMP are known to be excellent solvers.}. 	  
The solvers, however, are designed to detect infeasibility, and in the instances with infinite duality gap 
the dual is infeasible. 
Hence we would report that an instance is correctly solved, when Mosek or SDPA-GMP report 
dual infeasibility. However,  this did not happen for any of the instances.

We also tested the preprocessing method of \cite{permenter2014partial} and  Sieve-SDP 
\cite{zhu2019sieve} on the dual problems, then on the preprocessed problems we ran Mosek. 
Both these methods 
correctly preprocessed  all ``clean" instances, but could not preprocess the ``messy" instances. See the rows in Table \ref{table1} marked by 
PP+Mosek and Sieve-SDP+Mosek.  

\renewcommand*{\arraystretch}{1.3}
\begin{table}[ht]
\begin{center}
\begin{tabular}{ |c|c|c|c|c|}
\hline
& \multicolumn{2}{|c|}{Gap, single, finite} & \multicolumn{2}{c|}{Gap, single, infinite} \\
\cline{2-5}
& Clean & Messy & Clean & Messy   \\
\hline
Mosek & 1 & 1 & 0 & 0 \\ \hline
SDPA-GMP & 1 & 1 & 0 & 0  \\ \hline 
PP+Mosek  & 10 & 1 & 10 & 0 \\ \hline
Sieve-SDP + Mosek & 10& 1 & 10 & 0 \\
\hline
\end{tabular}
\end{center}
\caption{Computational results} 
\label{table1} 
\end{table}

We finally tested the exact SDP solver SPECTRA \cite{henrion2017spectra} on the
${\mathit gap}\_{\mathit single}\_{\mathit finite}\_{\mathit messy}\_2$  instance. SPECTRA cannot run on our SDPs as they are, since 
they do not satisfy 
an important algebraic assumption. However,  SPECTRA could compute and certify in exact arithmetic 
the optimal solution of the perturbed dual
\beq
\ba{rl}
\inf & B \bullet Y \\
s.t.    & A_i \bullet (Y  + \eps I) \, = \,\, c_i (i=1,2) \\
         & Y \succeq 0  
\ena
\eeq 
where $\epsilon$ was chosen as a small rational number. 

For example, with $\epsilon = 10^{-200}$ SPECTRA found and certified the optimal solution 
value as  $10.00000000$ in about two seconds of computing time.

\co{
We could possibly tackle our primal SDPs (but probably not the duals) by the facial reduction 
implementation in \cite{CheWolkSchurr:12}. It is not difficult to check that the singularity degree of the primal SDPs is just $1, \, $ and \cite{CheWolkSchurr:12} proves that their facial reduction algorithm can compute the minimal cone in a stable fashion. 
Indeed, \cite{Wolk:19} reported success on some of our instances in the $m=2$ case. 
}
The instances are stored in Sedumi format \cite{sturm1999using}, 
in which the roles of $c$ and $B$ are interchanged.
Each clean instance is given by 
\begin{itemize}
\item $A, \,$ which is a matrix with $m$ rows and $n^2$ columns, and the $i$th row of 
$A$ contains matrix $A_i$ of (\ref{p}) stretched out as a vector; 
\item $b, \,$ which is the $c$ in the primal (\ref{p}) , i.e., 
$b = 10 \cdot e_{n-1};$
\item $c, \,$ which is the right hand side $B$ of $(P),$ stretched out as a vector; 
\end{itemize}

These SDPs are available from the author's website.

\section{Conclusion}

We analyzed semidefinite programs with positive duality gaps, which, by common consent,
 is their 
most interesting and challenging pathology. 

We first dealt with the 
two variable case:  we   transformed two-variable SDPs into a standard form, that makes the 
positive gap  (if any) self-evident.  Second, we showed that the two variable case helps us understand 
  positive gaps in larger SDPs: the structure that causes a positive gap when $m = 2, \, $ often does the same in higher dimensions. 
  We  then investigated an intrinsic  parameter, the singularity degree of the duals of our SDPs, and proved that these are the largest that permit a positive gap. Finally, we  created a problem library of innocent looking, but very difficult SDPs, 
and showed that they are currently unsolvable by modern interior point methods.

Many interesting questions arise. 
First, and foremost, how do we solve the SDPs in our library?
It would be interesting to try  the algorithms of  \cite{liu2019new} or  \cite{hauenstein2018numerical} on the duals of our 
${\mathit gap}\_{\mathit single}\_{\mathit inf}$ SDPs, which are weakly infeasible.
We could also possibly adapt these algorithms to tackle  the ${\mathit gap}\_{\mathit single}\_{\mathit finite}$ 
instances (which have a finite duality gap). 

More generally, how do we solve SDPs with positive duality gaps? SDPs with infinity duality gap are known to arise in polynomial optimization, see for example  \cite{schweighofer2005optimization} and \cite{Waki:12}. 
Note that many SDPs with positive duality gaps may be ``in hiding:" when attempting to solve them,  solvers may just fail or report  an incorrect solution.
  
Second, how do we characterize positive gaps, when $m \geq 3$? 
Note that we have not completely understood even the $m=3$ case! 

Third, can we use the insight gained about positive gaps in SDPs to better understand positive gaps  
in other convex optimization problems? 

 We hope that this work will stimulate further research into these questions. 
 


	\appendix
	\section{Proofs of technical statements} 
	\label{section-technicalproofs} 
	
	\subsection{ Proof of Lemma \ref{rotate}} 
	\label{appendix-proof-lemma-rotate}
	
	Let $Q_1 \in \rad{r_1 \times r_1}$ be a matrix of orthonormal eigenvectors of $G_{11}, \, $ 
	$Q_2 \in \rad{r_2 \times r_2}$ a matrix of suitably normalized eigenvectors of $G_{22}, $ and $T_1 = Q_1 \oplus Q_2.$ 
	
	Then 
	$$
	T_1^T G T_1 = \left(    \ba{c|cc} 
	\Omega  & V & W \\ \hline 
	V^T        & I_s & 0 \\
	W^T       & 0 & 0 \ena \right), \, 
	$$
	with $\Omega  \in \sym{r_1}, \,$ $s$ is equal to the rank of $G_{22}, \,$ 
	and $V$ and $W$ are possibly nonzero.
	
	Next, let 
	$$
	T_2 = \left( \ba{c|cc} 
	I_{r_1} & 0 & 0 \\ \hline
	-V^T & I_{s}  & 0 \\
	0 & 0 & I_{r_2 - s}  \ena \right), \,\, {\rm then} \,\, 
	T_2^T T_1^T G T_1 T_2 = \left( \ba{c|cc} 
	\Omega  - VV^T  & 0 & W \\ \hline
	0  & I_s & 0 \\
	W & 0 & 0  \ena \right). 
	$$
	Finally, let $Q_3 \in \rad{r_1 \times r_1}$ be a matrix of orthonormal eigenvectors of $\Omega - VV^T, \,$ and $T_3 = Q_3 \oplus I_{r_2}, $ then 
	$T_3^T T_2^T T_1^T G T_1 T_2 T_3$ is in the required form.
	
	Also note that
	$$
	T_i^T \bpx I_{r_1} & 0 \\ 0 & 0 \epx T_i  = \bpx I_{r_1} & 0 \\ 0 & 0 \epx \,\, {\rm for} \,\, i=1,2,3, 
	$$
	hence \co{	 applying the rotation $T_i^T()T_i$ to $I_{r_1}  \oplus 0$ leaves this matrix the same 
	for all $i.$ Thus} 
	$$
	T := T_1 T_2 T_3
	$$
	will do.
	\qed
	
	\subsection{Proof of Claim \ref{Claim-D} } 
	\label{appendix-proof-Yi-Ai} 
		   
	 For brevity, let $L := \lin \, \{ A_1, \dots, A_m \}.$ We use induction. 
	We first prove (\ref{eqn-Yi-Ai}) for $i=1,$ so let 
	$Y_1 \in L \cap (\psd{n} \setminus \{ 0 \})$  and write $Y_1 = \sum_{j=1}^{m} \lambda_{j} A_j$ 
	with some $\lambda_{j}$ reals.  
	The last row of $Y_1$ is 
	$$
	(\lambda_{2}, \dots, \lambda_{n-1}, 0, 0).
	$$
	Since $Y_1 \succeq 0$ we deduce 
	$\lambda_{2} = \dots = \lambda_{n-1} = 0.$ Since $Y_1 \neq 0, \, $ 
	we deduce $\lambda_{1} > 0, \, $ so our claim
follows. 
	
	Next, suppose  (\ref{eqn-Yi-Ai}) holds for some $i$ such that $1 \leq i < m-1, \,$ and let 
	$(Y_1, \dots, Y_{i+1})$ be a strict facial reduction sequence whose members are all in $L.$ 
	We have 
	\beqast 
	Y_{i +1} & \in & (\psd{n} \cap Y_1^\perp \cap \dots \cap Y_i^\perp)^* \\
	             & =  & \bigl( 0 \oplus \psd{n-i} \bigr)^*,
	             \eeqast
	             where the equation follows, 
	             since 	$(Y_1, \dots, Y_{i})$  is also a strict facial reduction sequence, and using the induction hypothesis.
	             Thus the lower $(n-i) \times (n-i)$ block of 
	$Y_{i+1}$ is psd. Considering the last row of $Y_{i+1}$ and using a similar argument as before,
	we deduce that $Y_{i+1}$ is a linear 
	combination of $A_1, \dots, A_i, A_{i+1}$ only, i.e., 
	\beq \label{eqn-Yi+1-lambdaj} 
	Y_{i+1}  = \sum_{j=1}^{i+1} \lambda_{j} A_j
	\eeq
	with some $\lambda_{j}$ reals. 
	
	We claim that  $\lambda_{i+1} > 0. $ Indeed, $\lambda_{i+1} \geq 0$ since the lower right order $n-i$ 
	block of $Y_{i+1}$ is psd; and $\lambda_{i+1} =0$ would imply $Y_{i+1} \in ( 0 \oplus \psd{n-i})^\perp, \,$ 
	which would contradict  the assumption that 
	  	$(Y_1, \dots, Y_{i+1})$  is strict. Thus 
	  \co{	\beqast
	  	\psd{n} \cap Y_1^\perp \cap \dots \cap Y_i^\perp \cap Y_{i+1}^\perp & = & ( 0 \oplus \psd{n-i})  \cap Y_{i+1}^\perp \, {\rm (by \, the \, inductive \, hypothesis)} \\
																					  	                       & = & ( 0 \oplus \psd{n-i})  \cap (\lambda_{i+1} A_{i+1})^\perp  \\
																					  	                       & = & ( 0 \oplus \psd{n-i-1}), 
	  	           \eeqast 
	  	        } 
	  	           $$
	  	           \begin{array}{rcll} 
	  	           \psd{n} \cap Y_1^\perp \cap \dots \cap Y_i^\perp \cap Y_{i+1}^\perp & = & ( 0 \oplus \psd{n-i})  \cap Y_{i+1}^\perp &  {\rm (by \, the \, inductive \, hypothesis)} \\
	  	           & = & ( 0 \oplus \psd{n-i})  \cap (\lambda_{i+1} A_{i+1})^\perp  & ({\rm by } \,    \eqref{eqn-Yi+1-lambdaj})            \\ 
	  	           & = & ( 0 \oplus \psd{n-i-1}) & ({\rm by} \, \lambda_{i+1}>0), 
	  	           \end{array} 
	  	           $$
	  	so 	the proof is complete.

	
	\qed

	\subsection{ Proof of Lemma  \ref{lemma-double-ref}}    
	\label{subsection-proof-thm-doubleref} 
	For brevity, let 
	$$
	d := d \eqref{HD}  \, {\rm and} \, L := \lin \, \{ A_1, \dots, A_m \}.
	$$
	We first prove that 
	\beq \label{equation-Bmaxslack} 
	B \, {\rm is \, the \, maximum \, rank \,  psd \, matrix \, in \,} L. 
	\eeq
	To do so, assume that  
	$$
	B' \, = \, \sum_{j=1}^m \lambda_j A_j + \lambda B \succeq 0 
	$$
	has larger rank, where the $\lambda_i$ and $\lambda$ are suitable scalars. 
	Then $B'$ has a nonzero element outside its upper $r \times r$  block.
	Let  $\epsilon > 0$ be such that $|\lambda| \epsilon < 1, \,$ then 
	$$
	B'' := \dfrac{1}{1 + \lambda \epsilon} \bigl( B  + \epsilon B' \bigr) 
	$$
	is a slack in (\ref{p}) with rank larger than $r, \,$ a contradiction. Thus (\ref{equation-Bmaxslack}) follows.
	
	Next we prove that there is $A_1', \dots, A_{m}'$ such that \eqref{equation-frs-Ai-prop1} and \eqref{equation-frs-Ai-prop2} below hold: 
	\begin{eqnarray} \label{equation-frs-Ai-prop1} 
	(B, A_1', \dots, A_{m}') & \in & \fr(\psd{n}) \, {\rm \, is \, strict, \, and \, defines \, the \, minimal \, cone \, of \, } \eqref{HD},  \\
	\label{equation-frs-Ai-prop2} 
	A_1', \dots, A_m' & \in & L. 
	\end{eqnarray}
	These statements almost hold by definition: since $d = m+1, \, $ by definition, there is 
	\begin{eqnarray*} 
	(A_0', A_1', \dots, A_{m}') & \in & \fr(\psd{n}) {\rm \, which \, is \, strict, \, and \, defines \, the \, minimal \, cone \, of \, } \eqref{HD}, \, {\rm and} \\
	A_0', A_1', \dots, A_m' & \in & L + \lin \{ B \}.  
	\end{eqnarray*}
	By  (\ref{equation-Bmaxslack}) we have 
	$( \psd{n} \cap A_0'^{\perp})^* \subseteq ( \psd{n} \cap B^{\perp})^*,$   so $B$ is the best matrix to start a 
	facial reduction sequence with all members in $L + \lin \{ B \}. \, $ So we can replace $A_0'$ by $B, \, $ hence 
	 (\ref{equation-frs-Ai-prop1}) follows.  
	
	To ensure (\ref{equation-frs-Ai-prop2}) 
	let $i \in \{1, \dots, m \}$ and write $A_i' = \sum_{j=1}^{d-1} \lambda_j A_j + \lambda B$ for some 
	$\lambda_j$ and $\lambda$ reals. By definition, we have 
	\beq \label{equation-Aiprime} 
	A_i' \in (\psd{n} \cap B^\perp \cap A_1^{' \perp} \cap \dots \cap A_{i-1}^{' \perp})^*,
	\eeq
	thus a direct calculation shows 
		\beq \label{equation-Aiprime-2} 
		A_i' - \lambda B \in (\psd{n} \cap B^\perp \cap A_1^{' \perp} \cap \dots \cap A_{i-1}^{' \perp})^*. 
		\eeq
	thus subtracting $\lambda B$ from $A_i'$ maintains property (\ref{equation-Aiprime}).
	Doing this for $i=1,2 \dots, m$ keeps the sequence strict, and ensures that 
	(\ref{equation-frs-Ai-prop2}) holds. 
	
	Since  $(B, A_1', \dots, A_{m}')$ is strict, by Theorem 1 
	in \cite{liu2017exact}  they are 
	linearly independent. Thus we can reformulate (\ref{p}) using only 
	operations (\ref{exch}) and (\ref{trans}) in Definition \ref{definition-reform}
	to replace $A_i$ by $A_i'$ for $i=1, \dots, m, \,$ 
	where $A_1', \dots, A_{m}'$ are as in (\ref{equation-frs-Ai-prop1})
	(in the process we also replace the $c_i$ by suitable $c_i'.$ )
	
	Finally, by Lemma 2 in \cite{liu2017exact} there is an invertible 
	matrix $T$ of order $n$ such that
	$$
	(T^TB T, T^TA_1'T, \dots, T^TA_{m}'T) 
	$$
	is a regularized facial reduction sequence. We replace $B$ by $T^TBT$ and $A_i'$ by $T^T A_i'T$ for all $i$ 
	and this completes the proof. 

	\qed

	\nin{\bf Acknowledgement} I am grateful to Shu Lu and Minghui Liu for many helpful discussions and 
	to Quoc Tran Dinh and Yuzixuan Zhu for help in the computational experiments. 
	I also thank Simone Naldi for his generous help in using the SPECTRA package. 
	Special thanks are also due to Alex Touzov and Yuzixuan Zhu for a careful reading of the paper and to Alex Touzov 
	for his help with creating the graphics. 
	
	\bibliographystyle{plain}
\bibliography{mysdpMelody} 
\end{document}

%% file: imgs/pic_simpleEx1v2.tex
\usetikzlibrary{decorations.pathreplacing,angles,quotes}

\begin{tikzpicture}
\tikzstyle{every left delimiter}=[xshift=0.9ex]
\tikzstyle{every right delimiter}=[xshift=-0.9ex]
\newcommand{\ph}[1]{\phantom{#1}}

\matrix[
    matrix of math nodes,
    row sep=0.1ex,
    column sep=0.1ex,
    left delimiter=(,right delimiter=),
    nodes={text width=0.8em, text height=0.8em, text depth=0.5ex, align=center},
    every node/.style={scale=1.1}
    ] (m) at (-2.2,0)
    {
    1 & \ph{0} & \ph{0} \\
    \ph{0} & 0 & \ph{0} \\
    \ph{0} & \ph{0} & 0 \\
    };
    \begin{scope}[on background layer]
        \filldraw[red!30!orange] (m-1-1.north west) -- 
        (m-1-1.north east) -- (m-1-1.south east)-- (m-1-1.south west) --
        cycle;
    \end{scope} 

\matrix[
    matrix of math nodes,
    row sep=0.1ex,
    column sep=0.1ex,
    left delimiter=(,right delimiter=),
    nodes={text width=0.8em, text height=0.8em, text depth=0.5ex, align=center},
    every node/.style={scale=1.1}
    ] (m) at (0,0)
    {
    0 & \ph{0} & 1 \\
    \ph{0} & 1 & \ph{0} \\
    1 & \ph{0} & 0 \\
    };
    \begin{scope}[on background layer]
        \filldraw[cyan!70!blue] (m-1-1.north west) -- 
        (m-1-3.north east) -- (m-1-3.south east)-- (m-1-1.south east) -- (m-3-1.south east)--
        (m-3-1.south west) --
        cycle;
        \filldraw[red!30!orange] (m-2-2.north west) -- 
        (m-2-2.north east) -- (m-2-2.south east)-- (m-2-2.south west) --
        cycle;
    \end{scope} 

\matrix[
    matrix of math nodes,
    row sep=0.1ex,
    column sep=0.1ex,
    left delimiter=(,right delimiter=),
    nodes={text width=0.8em, text height=0.8em, text depth=0.5ex, align=center},
    every node/.style={scale=1.1}
    ] (m) at (2.2,0)
    {
    1 & \ph{0} & \ph{0} \\
    \ph{0} & 1 & \ph{0} \\
    \ph{0} & \ph{0} & 0 \\
    };
    \begin{scope}[on background layer]
        \foreach \x in {1,...,2}
        \filldraw[red!30!orange] (m-\x-\x.north west) -- 
        (m-\x-\x.north east) -- (m-\x-\x.south east)-- (m-\x-\x.south west) --
        cycle;
    \end{scope} 
    
    \draw[decoration={brace,mirror},decorate,thick]
        (-2.2-0.6,-0.8) -- node[below=2.5pt] {$A_1$} (-2.2+0.6,-0.8);
    
    \draw[decoration={brace,mirror},decorate,thick]
        (0-0.6,-0.8) -- node[below=2.5pt] {$A_2$} (0+0.6,-0.8);
    
    \draw[decoration={brace,mirror},decorate,thick]
        (2.2-0.6,-0.8) -- node[below=2.5pt] {$B$} (2.2+0.6,-0.8);
    
\end{tikzpicture}

%% file: imgs/pic_simpleEx23v4.tex
\begin{tikzpicture}
\tikzstyle{every left delimiter}=[xshift=0.9ex]
\tikzstyle{every right delimiter}=[xshift=-0.9ex]

\newcommand{\ph}[1]{\phantom{#1}}

\matrix[
    matrix of math nodes,
    row sep=0.1ex,
    column sep=0.1ex,
    left delimiter=(,right delimiter=),
    nodes={text width=0.8em, text height=0.8em, text depth=0.5ex, align=center},
    every node/.style={scale=0.95}
    ] (m) at (-2.6,0)
    {
    1 & \ph{0} & \ph{0} \\
    \ph{0} & 0 & \ph{0} \\
    \ph{0} & \ph{0} & 0 \\
    };
    \begin{scope}[on background layer]
        \filldraw[red!30!orange] (m-1-1.north west) -- 
        (m-1-1.north east) -- (m-1-1.south east)-- (m-1-1.south west) --
        cycle;
    \end{scope} 

\matrix[
    matrix of math nodes,
    row sep=0.1ex,
    column sep=0.1ex,
    left delimiter=(,right delimiter=),
    nodes={text width=0.8em, text height=0.8em, text depth=0.5ex, align=center},
    every node/.style={scale=0.95}
    ] (m) at (0,0)
    {
    0 & \ph{0} & 1 \\
    \ph{0} & 1 & \ph{0} \\
    1 & \ph{0} & 0 \\
    };
    \begin{scope}[on background layer]
        \filldraw[cyan!70!blue] (m-1-1.north west) -- 
        (m-1-3.north east) -- (m-1-3.south east)-- (m-1-1.south east) -- (m-3-1.south east)--
        (m-3-1.south west) --
        cycle;
        \filldraw[red!30!orange] (m-2-2.north west) -- 
        (m-2-2.north east) -- (m-2-2.south east)-- (m-2-2.south west) --
        cycle;
    \end{scope} 

\matrix[
    matrix of math nodes,
    row sep=0.1ex,
    column sep=0.1ex,
    left delimiter=(,right delimiter=),
    nodes={text width=0.8em, text height=0.8em, text depth=0.5ex, align=center},
    every node/.style={scale=0.95}
    ] (m) at (2.6,0)
    {
    1 & \ph{0} & \ph{0} \\
    \ph{0} & 1 & \ph{0} \\
    \ph{0} & \ph{0} & 0 \\
    };
    \begin{scope}[on background layer]
        \foreach \x in {1,...,2}
        \filldraw[red!30!orange] (m-\x-\x.north west) -- 
        (m-\x-\x.north east) -- (m-\x-\x.south east)-- (m-\x-\x.south west) --
        cycle;
    \end{scope} 

\matrix[
    matrix of math nodes,
    row sep=0.1ex,
    column sep=0.1ex,
    left delimiter=(,right delimiter=),
    nodes={text width=0.8em, text height=0.8em, text depth=0.5ex, align=center},
    every node/.style={scale=0.95}
    ] (m) at (-2.8-1.4,-2.1)
    {
    1 & \ph{0} & \ph{0} & \ph{0} \\
    \ph{0} & 0 & \ph{0} & \ph{0} \\
    \ph{0} & \ph{0} & 0 & \ph{0} \\
    \ph{0} & \ph{0} & \ph{0} & 0 \\
    };
    \begin{scope}[on background layer]
        \filldraw[red!30!orange] (m-1-1.north west) -- 
        (m-1-1.north east) -- (m-1-1.south east)-- (m-1-1.south west) --
        cycle;
    \end{scope} 

\matrix[
    matrix of math nodes,
    row sep=0.1ex,
    column sep=0.1ex,
    left delimiter=(,right delimiter=),
    nodes={text width=0.8em, text height=0.8em, text depth=0.5ex, align=center},
    every node/.style={scale=0.95}
    ] (m) at (-1.4,-2.1)
    {
    0 & \ph{0} & \ph{0} & 1 \\
    \ph{0} & 1 & \ph{0} & \ph{0} \\
    \ph{0} & \ph{0} & 0 & \ph{0} \\
    1 & \ph{0} & \ph{0} & 0 \\
    };
    \begin{scope}[on background layer]
        \filldraw[cyan!70!blue] (m-1-1.north west) -- 
        (m-1-4.north east) -- (m-1-4.south east)-- (m-1-1.south east) -- (m-4-1.south east)--
        (m-4-1.south west) --
        cycle;
        \filldraw[red!30!orange] (m-2-2.north west) -- 
        (m-2-2.north east) -- (m-2-2.south east)-- (m-2-2.south west) --
        cycle;
    \end{scope} 

\matrix[
    matrix of math nodes,
    row sep=0.1ex,
    column sep=0.1ex,
    left delimiter=(,right delimiter=),
    nodes={text width=0.8em, text height=0.8em, text depth=0.5ex, align=center},
    every node/.style={scale=0.95}
    ] (m) at (1.4,-2.1)
    {
    0 & \ph{0} & \ph{0} & \ph{0} \\
    \ph{0} & 0 & \ph{0} & 1 \\
    \ph{0} & \ph{0} & 1 & \ph{0} \\
    \ph{0} & 1 & \ph{0} & 0 \\
    };
    \begin{scope}[on background layer]
        \filldraw[cyan!70!blue] (m-1-1.north west) -- 
        (m-1-4.north east) -- (m-2-4.south east)-- (m-2-2.south east) -- (m-4-2.south east)--
        (m-4-1.south west) --
        cycle;
        \filldraw[red!30!orange] (m-3-3.north west) -- 
        (m-3-3.north east) -- (m-3-3.south east)-- (m-3-3.south west) --
        cycle;
    \end{scope}

\matrix[
    matrix of math nodes,
    row sep=0.1ex,
    column sep=0.1ex,
    left delimiter=(,right delimiter=),
    nodes={text width=0.8em, text height=0.8em, text depth=0.5ex, align=center},
    every node/.style={scale=0.95}
    ] (m) at (2.8+1.4,-2.1)
    {
    1 & \ph{0} & \ph{0} & \ph{0} \\
    \ph{0} & 1 & \ph{0} & \ph{0} \\
    \ph{0} & \ph{0} & 1 & \ph{0} \\
    \ph{0} & \ph{0} & \ph{0} & 0 \\
    };
    \begin{scope}[on background layer]
        \foreach \x in {1,...,3}
        \filldraw[red!30!orange] (m-\x-\x.north west) -- 
        (m-\x-\x.north east) -- (m-\x-\x.south east)-- (m-\x-\x.south west) --
        cycle;
    \end{scope}

\matrix[
    matrix of math nodes,
    row sep=0.1ex,
    column sep=0.1ex,
    left delimiter=(,right delimiter=),
    nodes={text width=0.8em, text height=0.8em, text depth=0.5ex, align=center},
    every node/.style={scale=0.95}
    ] (m) at (-3.2*2,-4.7)
    {
    1 & \ph{0} & \ph{0} & \ph{0} & \ph{0} \\
    \ph{0} & 0 & \ph{0} & \ph{0} & \ph{0} \\
    \ph{0} & \ph{0} & 0 & \ph{0} & \ph{0} \\
    \ph{0} & \ph{0} & \ph{0} & 0 & \ph{0} \\
    \ph{0} & \ph{0} & \ph{0} & \ph{0} & 0 \\
    };
    \begin{scope}[on background layer]
        \filldraw[red!30!orange] (m-1-1.north west) -- 
        (m-1-1.north east) -- (m-1-1.south east)-- (m-1-1.south west) --
        cycle;
    \end{scope}

\matrix[
    matrix of math nodes,
    row sep=0.1ex,
    column sep=0.1ex,
    left delimiter=(,right delimiter=),
    nodes={text width=0.8em, text height=0.8em, text depth=0.5ex, align=center},
    every node/.style={scale=0.95}
    ] (m) at (-3.2,-4.7)
    {
    0 & \ph{0} & \ph{0} & \ph{0} & 1 \\
    \ph{0} & 1 & \ph{0} & \ph{0} & \ph{0} \\
    \ph{0} & \ph{0} & 0 & \ph{0} & \ph{0} \\
    \ph{0} & \ph{0} & \ph{0} & 0 & \ph{0} \\
    1 & \ph{0} & \ph{0} & \ph{0} & 0 \\
    };
    \begin{scope}[on background layer]
        \filldraw[cyan!70!blue] (m-1-1.north west) -- 
        (m-1-5.north east) -- (m-1-5.south east)-- (m-1-1.south east) -- (m-5-1.south east)--
        (m-5-1.south west) --
        cycle;
        \filldraw[red!30!orange] (m-2-2.north west) -- 
        (m-2-2.north east) -- (m-2-2.south east)-- (m-2-2.south west) --
        cycle;
    \end{scope}

\matrix[
    matrix of math nodes,
    row sep=0.1ex,
    column sep=0.1ex,
    left delimiter=(,right delimiter=),
    nodes={text width=0.8em, text height=0.8em, text depth=0.5ex, align=center},
    every node/.style={scale=0.95}
    ] (m) at (0,-4.7)
    {
    0 & \ph{0} & \ph{0} & \ph{0} & \ph{0} \\
    \ph{0} & 0 & \ph{0} & \ph{0} & 1 \\
    \ph{0} & \ph{0} & 1 & \ph{0} & \ph{0} \\
    \ph{0} & \ph{0} & \ph{0} & 0 & \ph{0} \\
    \ph{0} & 1 & \ph{0} & \ph{0} & 0 \\
    };
    \begin{scope}[on background layer]
        \filldraw[cyan!70!blue] (m-1-1.north west) -- 
        (m-1-5.north east) -- (m-2-5.south east)-- (m-2-2.south east) -- (m-5-2.south east)--
        (m-5-1.south west) --
        cycle;
        \filldraw[red!30!orange] (m-3-3.north west) -- 
        (m-3-3.north east) -- (m-3-3.south east)-- (m-3-3.south west) --
        cycle;
    \end{scope}

\matrix[
    matrix of math nodes,
    row sep=0.1ex,
    column sep=0.1ex,
    left delimiter=(,right delimiter=),
    nodes={text width=0.8em, text height=0.8em, text depth=0.5ex, align=center},
    every node/.style={scale=0.95}
    ] (m) at (3.2,-4.7)
    {
    0 & \ph{0} & \ph{0} & \ph{0} & \ph{0} \\
    \ph{0} & 0 & \ph{0} & \ph{0} & \ph{0} \\
    \ph{0} & \ph{0} & 0 & \ph{0} & 1 \\
    \ph{0} & \ph{0} & \ph{0} & 1 & \ph{0} \\
    \ph{0} & \ph{0} & 1 & \ph{0} & 0 \\
    };
    \begin{scope}[on background layer]
        \filldraw[cyan!70!blue] (m-1-1.north west) -- 
        (m-1-5.north east) -- (m-3-5.south east)-- (m-3-3.south east) -- (m-5-3.south east)--
        (m-5-1.south west) --
        cycle;
        \filldraw[red!30!orange] (m-4-4.north west) -- 
        (m-4-4.north east) -- (m-4-4.south east)-- (m-4-4.south west) --
        cycle;
    \end{scope}

\matrix[
    matrix of math nodes,
    row sep=0.1ex,
    column sep=0.1ex,
    left delimiter=(,right delimiter=),
    nodes={text width=0.8em, text height=0.8em, text depth=0.5ex, align=center},
    every node/.style={scale=0.95}
    ] (m) at (3.2*2,-4.7)
    {
    1 & \ph{0} & \ph{0} & \ph{0} & \ph{0} \\
    \ph{0} & 1 & \ph{0} & \ph{0} & \ph{0} \\
    \ph{0} & \ph{0} & 1 & \ph{0} & \ph{0} \\
    \ph{0} & \ph{0} & \ph{0} & 1 & \ph{0} \\
    \ph{0} & \ph{0} & \ph{0} & \ph{0} & 0 \\
    };
    \begin{scope}[on background layer]
        \foreach \x in {1,...,4}
        \filldraw[red!30!orange] (m-\x-\x.north west) -- 
        (m-\x-\x.north east) -- (m-\x-\x.south east)-- (m-\x-\x.south west) --
        cycle;
    \end{scope}

\end{tikzpicture}

%% file: imgs/pic_simpleEx4v5.tex
\begin{tikzpicture}
\tikzstyle{every left delimiter}=[xshift=0.9ex]
\tikzstyle{every right delimiter}=[xshift=-0.9ex]

\newcommand{\ph}[1]{\phantom{#1}}

\matrix[
    matrix of math nodes,
    row sep=0.1ex,
    column sep=0.1ex,
    left delimiter=(,right delimiter=),
    nodes={text width=0.8em, text height=0.8em, text depth=0.5ex, align=center},
    every node/.style={scale=0.65}
    ] (m) at (-2.8,0)
    {
    1 & \ph{0} & \ph{0} & \ph{0} & \ph{0} \\
    \ph{0} & 0 & \ph{0} & \ph{0} & \ph{0} \\
    \ph{0} & \ph{0} & 1 & \ph{0} & \ph{0} \\
    \ph{0} & \ph{0} & \ph{0} & 0 & \ph{0} \\
    \ph{0} & \ph{0} & \ph{0} & \ph{0} & 0 \\
    };
    \begin{scope}[on background layer]
        \filldraw[red!30!orange] (m-1-1.north west) -- 
        (m-1-1.north east) -- (m-1-1.south east)-- (m-1-1.south west) --
        cycle;
        \filldraw[red!30!orange] (m-3-3.north west) -- 
        (m-3-3.north east) -- (m-3-3.south east)-- (m-3-3.south west) --
        cycle;
    \end{scope} 

\matrix[
    matrix of math nodes,
    row sep=0.1ex,
    column sep=0.1ex,
    left delimiter=(,right delimiter=),
    nodes={text width=0.8em, text height=0.8em, text depth=0.5ex, align=center},
    every node/.style={scale=0.65}
    ] (m) at (0,0)
    {
    0 & \ph{0} & \ph{0} & \ph{0} & 1 \\
    \ph{0} & 1 & \ph{0} & \ph{0} & \ph{0} \\
    \ph{0} & \ph{0} & 0 & \ph{0} & 1 \\
    \ph{0} & \ph{0} & \ph{0} & -1 & \ph{0} \\
    1 & \ph{0} & 1 & \ph{0} & 0 \\
    };
    \begin{scope}[on background layer]
        \filldraw[cyan!70!blue] (m-1-1.north west) -- (m-1-5.north east)--
        (m-1-5.south east) -- (m-2-3.north east)-- (m-2-3.south east) -- (m-3-5.north east)--
        (m-3-5.south east) -- (m-3-3.south east)--
        (m-5-3.south east) -- (m-5-3.south west)--
        (m-3-2.south east) -- (m-3-2.south west)--
        (m-5-1.south east) -- (m-5-1.south west)--
        cycle;
        \filldraw[red!30!orange] (m-2-2.north west) -- 
        (m-2-2.north east) -- (m-2-2.south east)-- (m-2-2.south west) --
        cycle;
        \filldraw[red!30!orange] (m-4-4.north west) -- 
        (m-4-4.north east) -- (m-4-4.south east)-- (m-4-4.south west) --
        cycle;
    \end{scope} 

\matrix[
    matrix of math nodes,
    row sep=0.1ex,
    column sep=0.1ex,
    left delimiter=(,right delimiter=),
    nodes={text width=0.8em, text height=0.8em, text depth=0.5ex, align=center},
    every node/.style={scale=0.65}
    ] (m) at (2.8,0)
    {
    1 & \ph{0} & \ph{0} & \ph{0} & \ph{0} \\
    \ph{0} & 1 & \ph{0} & \ph{0} & \ph{0} \\
    \ph{0} & \ph{0} & 1 & \ph{0} & \ph{0} \\
    \ph{0} & \ph{0} & \ph{0} & 0 & \ph{0} \\
    \ph{0} & \ph{0} & \ph{0} & \ph{0} & 0 \\
    };
    \begin{scope}[on background layer]
        \foreach \x in {1,...,3}
        \filldraw[red!30!orange] (m-\x-\x.north west) -- 
        (m-\x-\x.north east) -- (m-\x-\x.south east)-- (m-\x-\x.south west) --
        cycle;
    \end{scope} 

\matrix[
    matrix of math nodes,
    row sep=0.1ex,
    column sep=0.1ex,
    left delimiter=(,right delimiter=),
    nodes={text width=0.8em, text height=0.8em, text depth=0.5ex, align=center},
    every node/.style={scale=0.65}
    ] (m) at (-3.3*3/2,-2.3)
    {
    1 & \ph{0} & \ph{0} & \ph{0} & \ph{0} & \ph{0} & \ph{0} \\
    \ph{0} & 0 & \ph{0} & \ph{0} & \ph{0} & \ph{0} & \ph{0} \\
    \ph{0} & \ph{0} & 0 & \ph{0} & \ph{0} & \ph{0} & \ph{0} \\
    \ph{0} & \ph{0} & \ph{0} & 1 & \ph{0} & \ph{0} & \ph{0} \\
    \ph{0} & \ph{0} & \ph{0} & \ph{0} & 0 & \ph{0} & \ph{0} \\
    \ph{0} & \ph{0} & \ph{0} & \ph{0} & \ph{0} & 0 & \ph{0} \\
    \ph{0} & \ph{0} & \ph{0} & \ph{0} & \ph{0} & \ph{0} & 0 \\
    };
    \begin{scope}[on background layer]
        \filldraw[red!30!orange] (m-1-1.north west) -- 
        (m-1-1.north east) -- (m-1-1.south east)-- (m-1-1.south west) --
        cycle;
        \filldraw[red!30!orange] (m-4-4.north west) -- 
        (m-4-4.north east) -- (m-4-4.south east)-- (m-4-4.south west) --
        cycle;
    \end{scope} 

\matrix[
    matrix of math nodes,
    row sep=0.1ex,
    column sep=0.1ex,
    left delimiter=(,right delimiter=),
    nodes={text width=0.8em, text height=0.8em, text depth=0.5ex, align=center},
    every node/.style={scale=0.65}
    ] (m) at (-3.3*1/2,-2.3)
    {
    0 & \ph{0} & \ph{0} & \ph{0} & \ph{0} & \ph{0} & 1 \\
    \ph{0} & 1 & \ph{0} & \ph{0} & \ph{0} & \ph{0} & \ph{0} \\
    \ph{0} & \ph{0} & 0 & \ph{0} & \ph{0} & \ph{0} & \ph{0} \\
    \ph{0} & \ph{0} & \ph{0} & 0 & \ph{0} & \ph{0} & 1 \\
    \ph{0} & \ph{0} & \ph{0} & \ph{0} & 1 & \ph{0} & \ph{0} \\
    \ph{0} & \ph{0} & \ph{0} & \ph{0} & \ph{0} & 0 & \ph{0} \\
    1 & \ph{0} & \ph{0} & 1 & \ph{0} & \ph{0} & 0 \\
    };
    \begin{scope}[on background layer]
        \filldraw[cyan!70!blue] (m-1-1.north west) -- (m-1-7.north east)--
        (m-1-7.south east) -- (m-2-4.north east)-- (m-3-4.south east) -- (m-4-7.north east)--
        (m-4-7.south east) -- (m-4-4.south east)--
        (m-7-4.south east) -- (m-7-4.south west)--
        (m-4-3.south east) -- (m-4-2.south west)--
        (m-7-1.south east) -- (m-7-1.south west)--
        cycle;
        \filldraw[white] (m-2-2.north west) -- 
        (m-2-3.north east) -- (m-3-3.south east)-- (m-3-2.south west) --
        cycle;
        \filldraw[red!30!orange] (m-2-2.north west) -- 
        (m-2-2.north east) -- (m-2-2.south east)-- (m-2-2.south west) --
        cycle;
        \filldraw[red!30!orange] (m-5-5.north west) -- 
        (m-5-5.north east) -- (m-5-5.south east)-- (m-5-5.south west) --
        cycle;
    \end{scope} 

\matrix[
    matrix of math nodes,
    row sep=0.1ex,
    column sep=0.1ex,
    left delimiter=(,right delimiter=),
    nodes={text width=0.8em, text height=0.8em, text depth=0.5ex, align=center},
    every node/.style={scale=0.65}
    ] (m) at (3.3*1/2,-2.3)
    {
    0 & \ph{0} & \ph{0} & \ph{0} & \ph{0} & \ph{0} & \ph{0} \\
    \ph{0} & 0 & \ph{0} & \ph{0} & \ph{0} & \ph{0} & 1 \\
    \ph{0} & \ph{0} & 1 & \ph{0} & \ph{0} & \ph{0} & \ph{0} \\
    \ph{0} & \ph{0} & \ph{0} & 0 & \ph{0} & \ph{0} & \ph{0} \\
    \ph{0} & \ph{0} & \ph{0} & \ph{0} & 0 & \ph{0} & 1 \\
    \ph{0} & \ph{0} & \ph{0} & \ph{0} & \ph{0} & -1 & \ph{0} \\
    \ph{0} & 1 & \ph{0} & \ph{0} & 1 & \ph{0} & 0 \\
    };
    \begin{scope}[on background layer]
        \filldraw[cyan!70!blue] (m-1-1.north west) -- (m-1-7.north east)--
        (m-2-7.south east) -- (m-3-5.north east)-- (m-3-5.south east) -- (m-4-7.north east)--
        (m-5-7.south east) -- (m-5-5.south east)--
        (m-7-5.south east) -- (m-7-3.south east)--
        (m-5-3.south east) -- (m-5-3.south west)--
        (m-7-2.south east) -- (m-7-1.south west)--
        cycle;
        \filldraw[red!30!orange] (m-3-3.north west) -- 
        (m-3-3.north east) -- (m-3-3.south east)-- (m-3-3.south west) --
        cycle;
        \filldraw[red!30!orange] (m-6-6.north west) -- 
        (m-6-6.north east) -- (m-6-6.south east)-- (m-6-6.south west) --
        cycle;
    \end{scope}

\matrix[
    matrix of math nodes,
    row sep=0.1ex,
    column sep=0.1ex,
    left delimiter=(,right delimiter=),
    nodes={text width=0.8em, text height=0.8em, text depth=0.5ex, align=center},
    every node/.style={scale=0.65}
    ] (m) at (3.3*3/2,-2.3)
    {
    1 & \ph{0} & \ph{0} & \ph{0} & \ph{0} & \ph{0} & \ph{0} \\
    \ph{0} & 1 & \ph{0} & \ph{0} & \ph{0} & \ph{0} & \ph{0} \\
    \ph{0} & \ph{0} & 1 & \ph{0} & \ph{0} & \ph{0} & \ph{0} \\
    \ph{0} & \ph{0} & \ph{0} & 1 & \ph{0} & \ph{0} & \ph{0} \\
    \ph{0} & \ph{0} & \ph{0} & \ph{0} & 0 & \ph{0} & \ph{0} \\
    \ph{0} & \ph{0} & \ph{0} & \ph{0} & \ph{0} & 0 & \ph{0} \\
    \ph{0} & \ph{0} & \ph{0} & \ph{0} & \ph{0} & \ph{0} & 0 \\
    };
    \begin{scope}[on background layer]
        \foreach \x in {1,...,4}
        \filldraw[red!30!orange] (m-\x-\x.north west) -- 
        (m-\x-\x.north east) -- (m-\x-\x.south east)-- (m-\x-\x.south west) --
        cycle;
    \end{scope}

\matrix[
    matrix of math nodes,
    row sep=0.1ex,
    column sep=0.1ex,
    left delimiter=(,right delimiter=),
    nodes={text width=0.8em, text height=0.8em, text depth=0.5ex, align=center},
    every node/.style={scale=0.65}
    ] (m) at (-3.4*2,-5.4)
    {
    1 & \ph{0} & \ph{0} & \ph{0} & \ph{0} & \ph{0} & \ph{0} & \ph{0} & \ph{0} \\
    \ph{0} & 0 & \ph{0} & \ph{0} & \ph{0} & \ph{0} & \ph{0} & \ph{0} & \ph{0} \\
    \ph{0} & \ph{0} & 0 & \ph{0} & \ph{0} & \ph{0} & \ph{0} & \ph{0} & \ph{0} \\
    \ph{0} & \ph{0} & \ph{0} & 0 & \ph{0} & \ph{0} & \ph{0} & \ph{0} & \ph{0} \\
    \ph{0} & \ph{0} & \ph{0} & \ph{0} & 1 & \ph{0} & \ph{0} & \ph{0} & \ph{0} \\
    \ph{0} & \ph{0} & \ph{0} & \ph{0} & \ph{0} & 0 & \ph{0} & \ph{0} & \ph{0} \\
    \ph{0} & \ph{0} & \ph{0} & \ph{0} & \ph{0} & \ph{0} & 0 & \ph{0} & \ph{0} \\
    \ph{0} & \ph{0} & \ph{0} & \ph{0} & \ph{0} & \ph{0} & \ph{0} & 0 & \ph{0} \\
    \ph{0} & \ph{0} & \ph{0} & \ph{0} & \ph{0} & \ph{0} & \ph{0} & \ph{0} & 0 \\
    };
    \begin{scope}[on background layer]
        \filldraw[red!30!orange] (m-1-1.north west) -- 
        (m-1-1.north east) -- (m-1-1.south east)-- (m-1-1.south west) --
        cycle;
        \filldraw[red!30!orange] (m-5-5.north west) -- 
        (m-5-5.north east) -- (m-5-5.south east)-- (m-5-5.south west) --
        cycle;
    \end{scope}

\matrix[
    matrix of math nodes,
    row sep=0.1ex,
    column sep=0.1ex,
    left delimiter=(,right delimiter=),
    nodes={text width=0.8em, text height=0.8em, text depth=0.5ex, align=center},
    every node/.style={scale=0.65}
    ] (m) at (-3.4*1,-5.4)
    {
    0 & \ph{0} & \ph{0} & \ph{0} & \ph{0} & \ph{0} & \ph{0} & \ph{0} & 1 \\
    \ph{0} & 1 & \ph{0} & \ph{0} & \ph{0} & \ph{0} & \ph{0} & \ph{0} & \ph{0} \\
    \ph{0} & \ph{0} & 0 & \ph{0} & \ph{0} & \ph{0} & \ph{0} & \ph{0} & \ph{0} \\
    \ph{0} & \ph{0} & \ph{0} & 0 & \ph{0} & \ph{0} & \ph{0} & \ph{0} & \ph{0} \\
    \ph{0} & \ph{0} & \ph{0} & \ph{0} & 0 & \ph{0} & \ph{0} & \ph{0} & 1 \\
    \ph{0} & \ph{0} & \ph{0} & \ph{0} & \ph{0} & 1 & \ph{0} & \ph{0} & \ph{0} \\
    \ph{0} & \ph{0} & \ph{0} & \ph{0} & \ph{0} & \ph{0} & 0 & \ph{0} & \ph{0} \\
    \ph{0} & \ph{0} & \ph{0} & \ph{0} & \ph{0} & \ph{0} & \ph{0} & 0 & \ph{0} \\
    1 & \ph{0} & \ph{0} & \ph{0} & 1 & \ph{0} & \ph{0} & \ph{0} & 0 \\
    };
    \begin{scope}[on background layer]
        \filldraw[cyan!70!blue] (m-1-1.north west) -- (m-1-9.north east)--
        (m-1-9.south east) -- (m-2-5.north east)-- (m-4-5.south east) -- (m-5-9.north east)--
        (m-5-9.south east) -- (m-5-5.south east)--
        (m-9-5.south east) -- (m-9-5.south west)--
        (m-5-4.south east) -- (m-5-2.south west)--
        (m-9-1.south east) -- (m-9-1.south west)--
        cycle;
        \filldraw[white] (m-2-2.north west) -- 
        (m-2-4.north east) -- (m-4-4.south east)-- (m-4-2.south west) --
        cycle;
        \filldraw[red!30!orange] (m-2-2.north west) -- 
        (m-2-2.north east) -- (m-2-2.south east)-- (m-2-2.south west) --
        cycle;
        \filldraw[red!30!orange] (m-6-6.north west) -- 
        (m-6-6.north east) -- (m-6-6.south east)-- (m-6-6.south west) --
        cycle;
    \end{scope}

\matrix[
    matrix of math nodes,
    row sep=0.1ex,
    column sep=0.1ex,
    left delimiter=(,right delimiter=),
    nodes={text width=0.8em, text height=0.8em, text depth=0.5ex, align=center},
    every node/.style={scale=0.65}
    ] (m) at (-3.4*0,-5.4)
    {
    0 & \ph{0} & \ph{0} & \ph{0} & \ph{0} & \ph{0} & \ph{0} & \ph{0} & \ph{0} \\
    \ph{0} & 0 & \ph{0} & \ph{0} & \ph{0} & \ph{0} & \ph{0} & \ph{0} & 1 \\
    \ph{0} & \ph{0} & 1 & \ph{0} & \ph{0} & \ph{0} & \ph{0} & \ph{0} & \ph{0} \\
    \ph{0} & \ph{0} & \ph{0} & 0 & \ph{0} & \ph{0} & \ph{0} & \ph{0} & \ph{0} \\
    \ph{0} & \ph{0} & \ph{0} & \ph{0} & 0 & \ph{0} & \ph{0} & \ph{0} & \ph{0} \\
    \ph{0} & \ph{0} & \ph{0} & \ph{0} & \ph{0} & 0 & \ph{0} & \ph{0} & 1 \\
    \ph{0} & \ph{0} & \ph{0} & \ph{0} & \ph{0} & \ph{0} & 1 & \ph{0} & \ph{0} \\
    \ph{0} & \ph{0} & \ph{0} & \ph{0} & \ph{0} & \ph{0} & \ph{0} & 0 & \ph{0} \\
    \ph{0} & 1 & \ph{0} & \ph{0} & \ph{0} & 1 & \ph{0} & \ph{0} & 0 \\
    };
    \begin{scope}[on background layer]
        \filldraw[cyan!70!blue] (m-1-1.north west) -- (m-1-9.north east)--
        (m-2-9.south east) -- (m-3-6.north east)-- (m-4-6.south east) -- (m-5-9.north east)--
        (m-6-9.south east) -- (m-6-6.south east)--
        (m-9-6.south east) -- (m-9-5.south west)--
        (m-6-4.south east) -- (m-6-3.south west)--
        (m-9-2.south east) -- (m-9-1.south west)--
        cycle;
        \filldraw[white] (m-3-3.north west) -- 
        (m-3-4.north east) -- (m-4-4.south east)-- (m-4-3.south west) --
        cycle;
        \filldraw[red!30!orange] (m-3-3.north west) -- 
        (m-3-3.north east) -- (m-3-3.south east)-- (m-3-3.south west) --
        cycle;
        \filldraw[red!30!orange] (m-7-7.north west) -- 
        (m-7-7.north east) -- (m-7-7.south east)-- (m-7-7.south west) --
        cycle;
    \end{scope}

\matrix[
    matrix of math nodes,
    row sep=0.1ex,
    column sep=0.1ex,
    left delimiter=(,right delimiter=),
    nodes={text width=0.8em, text height=0.8em, text depth=0.5ex, align=center},
    every node/.style={scale=0.65}
    ] (m) at (3.4*1,-5.4)
    {
    0 & \ph{0} & \ph{0} & \ph{0} & \ph{0} & \ph{0} & \ph{0} & \ph{0} & \ph{0} \\
    \ph{0} & 0 & \ph{0} & \ph{0} & \ph{0} & \ph{0} & \ph{0} & \ph{0} & \ph{0} \\
    \ph{0} & \ph{0} & 0 & \ph{0} & \ph{0} & \ph{0} & \ph{0} & \ph{0} & 1 \\
    \ph{0} & \ph{0} & \ph{0} & 1 & \ph{0} & \ph{0} & \ph{0} & \ph{0} & \ph{0} \\
    \ph{0} & \ph{0} & \ph{0} & \ph{0} & 0 & \ph{0} & \ph{0} & \ph{0} & \ph{0} \\
    \ph{0} & \ph{0} & \ph{0} & \ph{0} & \ph{0} & 0 & \ph{0} & \ph{0} & \ph{0} \\
    \ph{0} & \ph{0} & \ph{0} & \ph{0} & \ph{0} & \ph{0} & 0 & \ph{0} & 1 \\
    \ph{0} & \ph{0} & \ph{0} & \ph{0} & \ph{0} & \ph{0} & \ph{0} & -1 & \ph{0} \\
    \ph{0} & \ph{0} & 1 & \ph{0} & \ph{0} & \ph{0} & 1 & \ph{0} & 0 \\
    };
    \begin{scope}[on background layer]
        \filldraw[cyan!70!blue] (m-1-1.north west) -- (m-1-9.north east)--
        (m-3-9.south east) -- (m-4-7.north east)-- (m-4-7.south east) -- (m-5-9.north east)--
        (m-7-9.south east) -- (m-7-7.south east)--
        (m-9-7.south east) -- (m-9-5.south west)--
        (m-7-4.south east) -- (m-7-4.south west)--
        (m-9-3.south east) -- (m-9-1.south west)--
        cycle;
        \filldraw[red!30!orange] (m-4-4.north west) -- 
        (m-4-4.north east) -- (m-4-4.south east)-- (m-4-4.south west) --
        cycle;
        \filldraw[red!30!orange] (m-8-8.north west) -- 
        (m-8-8.north east) -- (m-8-8.south east)-- (m-8-8.south west) --
        cycle;
    \end{scope}

\matrix[
    matrix of math nodes,
    row sep=0.1ex,
    column sep=0.1ex,
    left delimiter=(,right delimiter=),
    nodes={text width=0.8em, text height=0.8em, text depth=0.5ex, align=center},
    every node/.style={scale=0.65}
    ] (m) at (3.4*2,-5.4)
    {
    1 &  &  &  &  &  &  &  &  \\
     & 1 &  &  &  &  &  &  &  \\
     &  & 1 &  &  &  &  &  &  \\
     &  &  & 1 &  &  &  &  &  \\
     &  &  &  & 1 &  &  &  &  \\
     &  &  &  &  & 0 &  &  &  \\
     &  &  &  &  &  & 0 &  &  \\
     &  &  &  &  &  &  & 0 &  \\
     &  &  &  &  &  &  &  & 0 \\
    };
    \begin{scope}[on background layer]
        \foreach \x in {1,...,5}
        \filldraw[red!30!orange] (m-\x-\x.north west) -- 
        (m-\x-\x.north east) -- (m-\x-\x.south east)-- (m-\x-\x.south west) --
        cycle;
    \end{scope}

\end{tikzpicture}